\documentclass[cup7b]{cupbook}
\input{defs.ttt}

\begin{document}

\chapter{Generalised Proof-Nets for Compact Categories with
  Biproducts}
\author{Ross Duncan}
\begin{abstract}
  Just as conventional functional programs may be understood as proofs in
  an intuitionistic logic, so quantum processes can also be viewed as
  proofs in a suitable logic.  We describe such a logic, the logic of
  compact closed categories and biproducts, presented both as a
  sequent calculus and as a system of proof-nets.  
  This logic captures much of the necessary structure needed to
  represent quantum processes under classical control, while remaining
  agnostic to the fine details.  We demonstrate how to represent
  quantum processes  
  as proof-nets, and show that the dynamic behaviour of a quantum
  process is captured by the cut-elimination procedure for the logic.
  We show that the cut elimination procedure is strongly
  normalising: that is, that every legal way of simplifying a proof-net
  leads to the same, unique, normal form.  Finally, taking some initial
  set of operations as 
  non-logical axioms, we show that that the resulting category of
  proof-nets  is a representation of the free compact 
  closed category with biproducts generated by those operations.
\end{abstract}


\section{Introduction}
\label{sec:introduction}

\subsection{Logic, Processes, and Categories}
\label{sec:categ-proofs-proc}

Birkhoff and von Neumann initiated the logical study of quantum 
mechanics in their 1936 paper \cite{Birkhoff1936The-Logic-of-Qu}.
They constructed a logic by assigning a proposition letter to each
observable property of a given quantum system, and studied 
negations, conjunctions, and disjunctions of these properties.  The resulting 
lattice is non-distributive, and so the heart of what is called
``quantum logic''is the study of various kinds of non-distributive
lattices.
These traditional quantum logics suffer from a number of defects.
Firstly, they are monolithic:  there is no way to derive 
the properties of a composite system from the properties of its
parts.  Each system has its own associated lattice which can only rarely
be related to those of other systems.  Further, these systems are seen
statically:  a system that undergoes some dynamical change is
 a new system.  Finally, the failure of compositionality is connected
 with the fact 
that quantum logic has no decent notion of implication
\cite{Smets2001The-Logic-of-Ph};  hence we have 
a logic which has a notion of validity, but no concept of
\emph{inference} or \emph{proof}.

These limitations form a serious obstacle to
the use of Birkhoff-von Neumann-style quantum logic to study
\emph{interacting} quantum 
systems.  If by ``system'' we understand the ubiquitous
qubit, then it is precisely such interactions which form 
the basis of quantum informatics.  Indeed, the principal 
concern of the computer scientist\footnote{At least: the principal
  concern of \emph{many} computer scientists.}---how to soundly
construct large 
systems from smaller ones---is exactly where quantum logic is
weakest.  

In this article I will describe a very different kind
of logic which can address these questions.  This logic is called
\emph{tensor-sum logic} and it shares many features with linear logic
\cite{Girard:LinLog:1987}, which has been widely studied in computer
science and structural proof theory.

We proceed in accordance with an old tradition in computer science,
that of linking 
computational systems and logics, of claiming that a certain logic ``is
the same as'' some formal computing machine.    The archetype of this
approach is the Curry-Howard correspondence between intuitionistic
natural deduction
\cite{Gentzen1935Untersuchungen,Prawitz1965Natural-Deducti} and the 
simply-typed $\lambda
$-calculus \cite{Curry1958Combinatory-Log,Howard:1980kb}.  The basic 
insight is that the inference 
rules of the logic are essentially the same as the constructors used
to form a $\lambda$-term.  The valid types of the $\lambda$-calculus are
nothing more than the theorems of the logic, and more importantly
every proof represents a $\lambda$-term. We can thus view the logic
itself as a computational system, with the  important proviso that the
objects of interest are the \emph{proofs} and not the theorems.
This correspondence is not just skin deep.  Recall that the $\beta$-reduction
relation between $\lambda $-terms expresses the \emph{execution} behaviour of
the calculus; we consider terms to be computationally equivalent when
they reduce to the same $\beta$-normal form.   This dynamical aspect
of the $\lambda $-calculus corresponds exactly to the normalisation procedure
for natural deduction proofs or, in the sequent calculus setting, to
the \emph{cut-elimination} procedure
\cite{Gentzen1935Untersuchungen}.  In our presentation of tensor-sum
logic we will take this correspondence as given, and treat
the proof theoretic presentation as a computational system in its own
right.  The normalisation procedure gives the computational dynamics.


There is a further correspondence that we engage with rather
more
seriously: that between the proofs in a given logic and the arrows in
a particular class of categories
\cite{Lambek1968Deductive-Syste,Lambek1969Deductive-Syste}.  We can equate
the propositions of a logic $L$ with the objects of 
a category $\catC$;  for each proof of, say, a proposition $B$ from a
premise $A$, we define a corresponding arrow $f : A \to B$ in the
category.
Stated so blandly, we have little structure to work 
with, so we go further, and demand that the logical connectives
are represented by functors on $\catC$. The natural
transformations between these functors then give rise to the inference
rules of the logic.   In this way an arrow in \catC may be constructed
for each proof in $L$;  we say that \catC is  a \emph{denotational
  model} of $L$ if, whenever two proofs  share the same normal
form, they have equal denotations in \catC.   The categorical model
gives an extensional account of the computational dynamics represented
in the term language.  The general schema of this tripartite relation
is shown in Table~\ref{tab:CHLelements}.

\begin{table}[htb]
  \centering 
  \begin{tabular}{||c|c|c||} \hline
    \textbf{Computation} & \textbf{Logic} & \textbf{Category} \\ \hline
    types & formulae & objects \\
    terms  & proofs & arrows \\
    type formers & connectives & functors \\
    term constructors & inference rules & natural transformations \\ \hline
  \end{tabular}
  \caption{Curry-Howard-Lambek correspondences}
  \label{tab:CHLelements}
\end{table}

To return to our earlier example, the simply-typed
$\lambda $-calculus forms such a triple with intuitionistic natural
deduction and the class of Cartesian closed
categories \cite{Lambek1986Introduction-to}.  Another example is
provided by intuitionistic multiplicative linear logic
\cite{Girard:LinLog:1987}; this logical system corresponds to the
\emph{linear} $\lambda$-calculus  and to the class of *-autonomous
categories \cite{Barr:*-AutoCats:1979,Barr:*-AutoLL:1991}.  The
pattern is quite general, and more examples can easily be found.

This relationship between a logic and its categorical models can be
made exact.  In the above description we embedded
the connectives and inference rules of $L$ into \catC using only
functors and natural transformations: the logical structure is
agnostic with respect to  the concrete elements of the category.
Hence any category with the requisite functors can provide a model, up
to some assignment of the basic proposition letters.  Of course such a
category may well contain other objects or arrows which do not
correspond to anything in the logic we are trying to model.  To make
the correspondence exact we must be able to translate from the
\emph{free category}  (with appropriate structure) faithfully back to
our logic.  That is, we must find an injective translation from the
arrows of \catC onto the cut-free proofs of $L$.
In the case of our generalised proof-nets, that is indeed possible,
but for a simpler example, consider the simply-typed
$\lambda $-calculus with just one ground type; then the category of
its terms is the free Cartesian closed category  generated by the
category $\mathbf{1}$, which has only one object, and one identity
arrow.

The choice of generator can be rather important.  By choosing
a discrete category (i.e. a set) the resulting logic will have that set
as its propositional variables.  By choosing a category with
non-trivial arrows, we introduce non-logical axioms:  if cut-elimination for
the resulting logic is to be retained we must lift the composition
operation of \catC into the cut-elimination procedure of $L$ itself.  If other
equations between arrows are required these too must be hoisted
into the logic.  In the case we consider
here, the situation is even worse: we will take as a generator a
\emph{compact symmetric polycategory} \cite{Duncan:thesis:2006}.  This
 esoteric creature will be described in a later section, but 
for now we note that the intricacy of the required composition
forces the adoption of the proof-net formalism, an illustration of the
power of graphical methods over conventional syntax.

The development in the following sections will be the reverse
of the exposition above.  We first describe the
mathematical basis of quantum computation in its concrete
setting---finite dimensional Hilbert spaces---then we identify certain
structures of the category \fdHilb which are the
essential features for carrying out quantum computation.   
Next we
present the syntax of tensor-sum logic, and finally prove that the
category of generalised proof-nets is the free compact closed
category with biproducts generated by a compact symmetric
polycategory.  

\paragraph{Prior Work}
\label{sec:related-work}
The original formulation of quantum mechanics in terms of compact
closed categories  with biproducts was due to Abramsky and Coecke
\cite{AbrCoe:CatSemQuant:2004}.  An early attempt to formalise quantum
computations in terms of proof-nets for \MLL was
\cite{Duncan:BellStateMLL:2004}.  The first description of 
 a logic based on compact closed categories and biproducts was given
 by Abramsky and the author in \cite{AbrDun:CQLv2:2004},  however this
 logic is  essentially restricted to quantum systems with only bipartite
 entanglement.  This restriction was lifted in
 \cite{Duncan:thesis:2006} via the use of polycategories, although
 only the multiplicative fragment of the logic was treated.
The present article is essentially a fusion of the previous two:
giving a complete presentation of two-sided proof-nets with
generalised axioms, and both the tensor and sum connectives.
A notable distinction between the treatment of the biproduct here as
compared to that of \cite{AbrCoe:CatSemQuant:2004} is that here the
biproduct is freely generated, hence the our treatment is closer to
that of \cite{Selinger:TowardQPL:2003};  we will discuss the vexed
position of this connective at the  end of Section
\ref{sec:quantum-computation}.

\subsection{Quantum Mechanics Concretely, and Abstractly}
\label{sec:quantum-computation}

The main work of this article is to characterise the structure of
certain kinds of category in terms of proof-nets.  The particular
categories we are interested in are\emph{ compact closed categories with
biproducts}, and in the next section we'll go into considerable detail
defining, characterising, and giving some of the basic properties of
such categories.   In this section we describe, at a more intuitive
level, how the structure relates to quantum mechanics and quantum
computation  in particular.

We start with a schematic description of quantum mechanics.  Since we
are interested in quantum computation, we will restrict 
ourselves to quantum systems with finite dimensional state spaces.  Consider
the following axioms.
\begin{enumerate}
\item To each quantum system we associate a finite dimensional Hilbert
  space, its 
  \emph{state space};  the possible states of the system are unit
  vectors in the state space, modulo a global phase factor.
\item If two systems are combined, then their joint state space is the
  \emph{tensor product} of the two state spaces.
\item Measuring a quantum system is \emph{non-deterministic};
  possible outcomes are the eigenvectors of some self-adjoint
  operator on the state space, and the probability of observing a
  particular outcome depends on the inner product of the current state
  and the eigenvector for that outcome.  After the measurement the
  state is updated to match the observed vector\footnote{We have
    implicitly assumed that the measurement is non-degenerate.},
  assuming it is not destroyed by the act of measuring.
\item For any discrete time step, the next state of the system is
  determined by a \emph{unitary map} on its state space.
\end{enumerate}
Of course, we have omitted some important details, but the above axioms
approximate the level of abstraction of our categorical
formalisation.

The first thing we should note that the \fdhilb, the  category of finite 
dimensional  Hilbert spaces and linear maps, is the ``natural''
category in which to formalise the above axioms\footnote{As the first
  and last axioms suggest \fdhilb is actually too big:  it contains
  vectors and maps which do not correspond to anything in quantum
  mechanics.  The general program of describing quantum mechanics in
  categorical terms aims to find the \emph{minimal} structure
  required.}.  The second point is that \fdhilb is compact closed  and
has biproducts. Thirdly, the above axioms can be rephrased in terms of
the categorical structure alone.

Before showing how the ingredients of the axioms translate into
categorical terms, let's dispose of some unnecessary baggage. Consider a
linear map $\psi : \mathbb{C} \to \mathcal{H}$.  Since it is linear,
and its domain is one-dimensional, its value is fixed by its value
on $1$, hence maps of this type are in 1-1 correspondence with vectors
of $\mathcal{H}$.  Hence we can forget about vectors and talk only of
linear maps.

Now consider measurements.  There are three parts to a quantum
measurement:  the non-deterministic possibilities, the calculation of
the probabilities, and the updated state.  We will deal with these in
reverse order.  The new state of a measured system depends only on
which outcome $i$ of the measurement happened, hence it is just a new
state $\ket{\phi_i}$, with no particular relation to the old one.  Of
course, the measurement process may destroy the system, in which case
there is no new state.  To calculate the probability of seeing outcome
$i$ when we are in state $\ket{\psi}$ we must calculate the inner
product $\innp{\phi_i}{\psi}$.  This process too can be seen as a
linear map, namely  the projection map $\bra{\phi_i}: \mathcal{H} \to
\mathbb{C}$; when composed with $\ket{\psi}$ this yields the inner
product.  Hence the state transformation associated to the $i$th
outcome of some measurement is described by the map:
\[
\mathcal{H} \rTo^{\bra{\phi_i}} 
\mathbb{C}  \rTo^{\ket{\phi_i}} 
\mathcal{H}\;.
\]
Note that this is desired transformation even in the  case when only
part of a composite system is measured.  For the purposes of this
work, calculating the probabilities is not so important, but the
transformation of the state is essential.  

Finally let's consider the non-deterministic aspect of measurement.
The main point is that there several possible outcomes \emph{and we
  know which one happened}.   Suppose we perform a measurement with
two outcomes; we can view this process as map of type
\[
\mathcal{H} \to \mathcal{H} \oplus \mathcal{H}
\]
where the two sides of the direct sum correspond to the two ``possible
worlds'' induced by the two outcomes of the measurement.  One can then
represent conditional operations by acting on only one subspace or the
other; for example $f\oplus g$ behaves as the linear map
$f:\mathcal{H} \to \mathcal{H}$ if the first outcome was observed, and
$g:\mathcal{H} \to \mathcal{H}$ if the second was observed.  

How then can we write the axioms of quantum mechanics in the language
of compact closed categories and biproducts?  Firstly, compact closed
categories are equipped with a tensor product, and this tensor product
has a neutral element $I$.  In \fdhilb the tensor product is simply
the usual Kronecker product, and its neutral element is the base
field, namely $\mathbb{C}$.  This gives us the first two axioms:
\begin{itemize}
\item To each quantum system we associate an object $A$, its state
  space;  the possible states of the system are given by arrows $\psi
  : I \to A$.
\item If two systems with state spaces $A$ and $B$ respectively are
  combined, then their joint state space is $A\otimes B$.
\end{itemize}
In \fdhilb the biproduct is the direct sum of Hilbert spaces, and this
will allow the formalisation of measurements.
\begin{itemize}
\item An $n$-outcome measurement of a quantum system whose state space
  is $A$ is represented by an arrow
  \[
  m: A \rTo \bigoplus_i B_i 
  \]
  where each of the the projections $\pi_i\circ m : A \to B_i$ factors
  as
  \[
  A \rTo I \rTo B_i\;.
  \]
\end{itemize}
This is very general notion of measurement. The common cases are when
$B_i = I$ and the original system is destroyed; and, when $B_i = A$
when the original system is preserved.

Since we are interested in formalising as much of quantum mechanics as
possible within one category we will not restrict state
transformations to unitary evolutions;  note that a measurement is a
valid transformation of a quantum state which is  not unitary.  Hence
the last axiom is simply the following.
\begin{itemize}
\item A quantum system $A$ may transform to another system $B$ by
  means  of any arrow $f:A\to B$.
\end{itemize}
\noindent
We have not yet mentioned the role of the compact structure of the
category.  While not required to paraphrase the axioms, the compact structure
 plays an important role in capturing quantum phenomena.  Recall that
 the tensor product of two vector spaces contains points, $\Psi : I
 \to A \otimes B$ which cannot be factored into a pair of vectors
$\psi_1 : I \to A$ and $\psi_2 : I \to B$.  The such quantum states
are called \emph{entangled} and they a central role in quantum
computation.  The compact structure guarantees the existence of
certain entangled states, namely \emph{Bell states} for every finite
dimensional Hilbert space:
\[
\eta_A : I \to A^* \otimes A\;.
\]
If $A$ is two-dimensional, i.e. a qubit, then the corresponding vector is
\[
\eta_A(1) = \ket{00} + \ket{11}\;.
\]
Further more, the compact structure also provides a projection onto
this state
\[
\epsilon_A : A \otimes A^* \to I
\]
hence we can define measurements on Bell states.  These two operations
will allow many more entangled states to be defined.

This completes our impressionistic description of how quantum
mechanics may be formalised in the categorical setting.
 Before moving on, it worth pointing out what has been
excluded from our formalisation.  Perhaps the most striking omission
in moving between the concrete axioms and the abstract is the concept
of \emph{unitarity}.   

The abstract formulation of quantum mechanics described here is
derived from that introduced in
\cite{AbrCoe:CatSemQuant:2004} which uses 
\emph{strongly} compact closed  categories.  Also called $\dag$-compact,
these categories are equipped with a contravariant involutive functor
which sends each map $f$ to its adjoint $f^\dag$, and has no effect on
objects.  This functor can then be used to define unitarity and the
inner product.  In this article, we focus on freely
constructing the compact closed and biproduct structure from some
underlying category of generators.  One could 
consider the case when these structures cohere with the $(\cdot)^\dag$
operation, giving a $\dag$-compact category with $\dag$-biproducts;
however the only difference here between the $\dag$-structure and the
original is that  the structural isomorphisms are required to be unitary.  That
is to say that the only \emph{new} maps which are introduced are the
adjoints of the generators.  Hence we can simply enlarge the class of
generators beforehand, and thereafter ignore the $\dag$-structure.
Of course, when working on concrete examples it is important to be
aware of the adjoints, and the equational theory of the generators
more generally, but that is not the focus of the present article.

The other important deviation from usual quantum mechanics is that we
have been extremely liberal about measurements.  In 
particular we do not make any restriction on the number of outcomes a
measurement may have.  Of course, in quantum mechanics the outcomes
are the spectrum of some operator, and hence are bounded by the
dimension of the space.  Considering these issues would take us too
far afield but \cite{AbrCoe:CatSemQuant:2004} has one approach;
a more recent categorical treatment of quantum observables is found in
\cite{Coecke:2008nx}.  In any case, it seems unlikely that the
structure of quantum measurements---being fundamentally connected to
the bases of the underlying space---will yield to a description in terms
of natural transformations of some functors.

\subsubsection*{Some Remarks on the Biproduct}
\label{sec:some-remarks-bipr}

In their original paper \cite{AbrCoe:CatSemQuant:2004} Abramsky and
Coecke used the biproduct of \fdHilb in two roles:  firstly, to encode
classical branching, as described above; and secondly, to construct
\emph{bases} for the underlying space.  In particular, they define
state preparations and destructive measurements as isomorphisms of the forms
\[
\text{base} : \bigoplus_i I \to A
\quad\text{and}\quad
\text{meas} : A \to I \bigoplus_i\;.
\]
This second use of the biproduct has been criticised by later works
\cite{Coecke2005De-linearizing-,Selinger:dagger:2005} on two main
accounts.  In the original approach the composite
\[
A \rTo^{\text{ meas }} \bigoplus_i I \rTo^{\text{ meas }} A
\]
yields the identity map, contradicting physical reality---a real
experiment would transform a pure state to a mixed state, something
not handled within this simple framework.  More importantly, when
moving from a ``vector space'' setting like \fdHilb to a
``projective'' setting, such as 
Selinger's CPM construction or Coecke's WProj, the direct sum of the
underlying space no longer yields a biproduct.  The only option is to
construct the biproduct as formal vectors and matrices.  The
works cited above show that, in the projective setting, if there is a
biproduct then the scalars are essentially restricted to \emph{probabilities}
rather than \emph{amplitudes}.  The immediate consequence is that we
must give up any hope of using the additive structure to encode
interference effects:  we are essentially restricted to a classical
probabilistic setting.

The approach to biproducts taken in this work is absolutely consonant
with these restrictions.  We construct both the multiplicative and
additive structures freely, and hence the scalars are simply a (free)
semiring.  The theory of processes thus produced is much like that
introduced in \cite{Selinger:TowardQPL:2003}, based on classically
controlled quantum operations.

\subsection{An Example Proof-net}
\label{sec:an-example-proof}

Sections \ref{sec:tensor-sum-logic} and \ref{sec:gener-proof-nets}
will introduce tensor-sum logic, and its proof-net notation.  Since
those sections will focus on the technical details of the
formalism we present now an illustrative example of how proof-nets can
be used to to model quantum processes.

We will describe an old favourite:  the quantum teleportation
protocol \cite{BBCJW:1993:teleport}.  The sketch of the  protocol goes
like this:
 Two parties, Alice and Bob, initially share an entangled pair of
  particles in a Bell state,
  \[
  \ket{\text{Bell}} =  \frac{\ket{00} + \ket{11}}{\sqrt{2}}
  \]
The parties then separate, and at some later point Alice wants to send a
qubit to Bob, but unfortunately she has only a classical channel.
However, it is still possible to transmit the qubit by using the shared
entanglement between the two parties.

To proceed, Alice performs a joint measurement on the qubit she wishes
to transmit together with her half of the  entangled pair.  She
measures in the Bell basis, so her state will be projected onto one of
the following vectors:
\begin{gather*}
  \ket{\text{Bell}_1} =  \frac{\ket{00} + \ket{11}}{\sqrt{2}}
  \quad
  \ket{\text{Bell}_Z} =  \frac{\ket{00} - \ket{11}}{\sqrt{2}}
\\
  \ket{\text{Bell}_X} =  \frac{\ket{01} + \ket{10}}{\sqrt{2}}
  \quad
  \ket{\text{Bell}_Y} =  \frac{\ket{01} - \ket{10}}{\sqrt{2}}
\end{gather*}
These states are all entangled, and further, each of them can be
produced by starting with $\ket{Bell_1}$  (aka $\eta$) and applying
one of the Pauli operators;  hence  we can associate a
Pauli operator to each outcome of the measurement.  In order to
complete the protocol, Alice transmits a classical message to Bob,
saying which of the  four outcomes she observed.  Bob then applies the
corresponding Pauli operation to his qubit and---as if by magic---it
is now in the state that Alice wished to transmit.

We now show how this protocol can be represented in term of proof-nets
using the compact closed  structure and biproducts.  By normalising the
proof-net we will effectively simulate the execution of the  protocol.

We start with a \emph{premise} representing Alice's input, and a unit
link, representing the initial shared Bell state
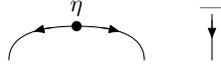
\[
  \raisebox{-0.5\height}{\scalebox{0.9}{
      \includegraphics{images/newproofnets_36.mps}}}
\]
Note that in these diagrams time flows from the top to the bottom:
input at the top of the page, outputs at the bottom.  The right two
qubits are taken to belong to Alice, the leftmost belonging to Bob.
The next element is the Bell basis measurement.  We will assume this is a
destructive measurement, so the four possible transformations are
simply projections.  To indicate that these form an exclusive
choice, we put them in a box, as shown below.
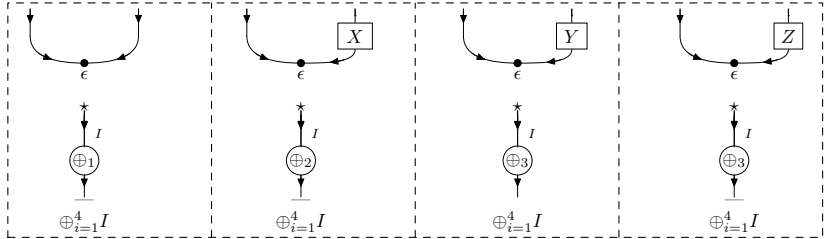
\[
  \raisebox{-0.5\height}{\scalebox{0.9}{
      \includegraphics{images/newproofnets_37.mps}}}
\]
Notice the output of type $\oplus_i I$ serves simply to indicate which
outcome occurred

Finally we consider Bob's correction.  Since his behaviour is
conditional on a classical input, he has a box with an input of type 
$\oplus_i I$ as shown below.
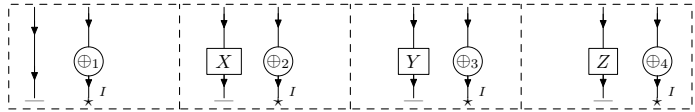
\[
  \raisebox{-0.5\height}{\scalebox{0.9}{
      \includegraphics{images/newproofnets_38.mps}}}
\]
Putting it all together we have the following picture:
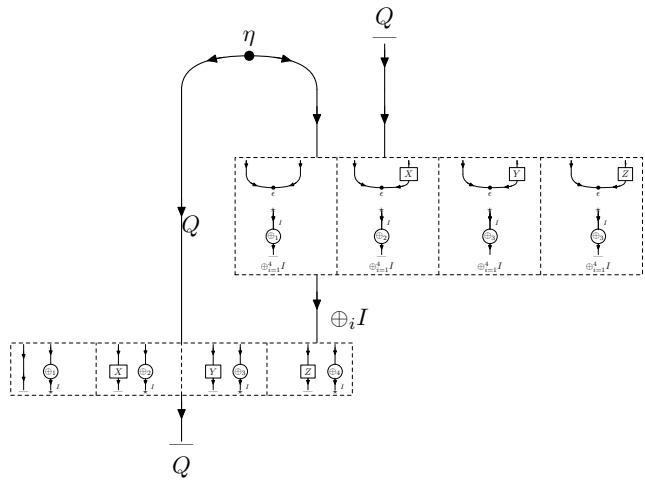
\[
  \raisebox{-0.5\height}{\scalebox{0.9}{
      \includegraphics{images/newproofnets_39.mps}}}
\]
Now we can begin to simulate the protocol.  The first step is to
resolve the non-determinism of Alice's measurement.  We do this by
``opening the box'', essentially making four copies of the whole
system, one for each possible outcome of the measurement.
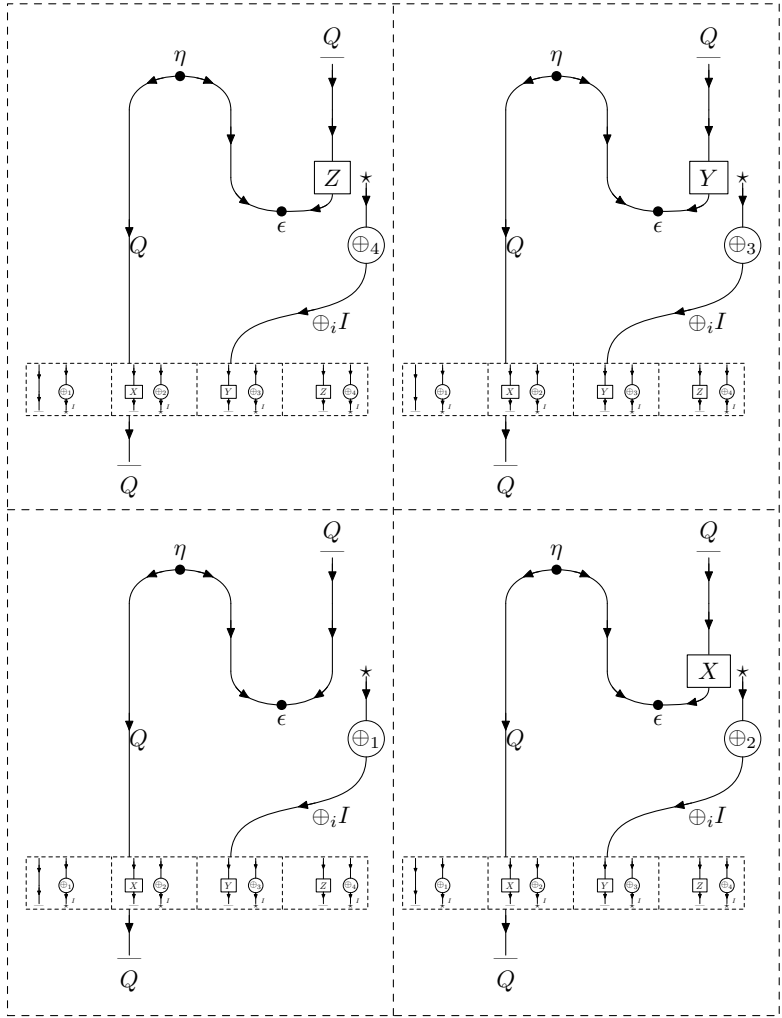
\[
  \raisebox{-0.5\height}{\scalebox{0.9}{
      \includegraphics{images/newproofnets_40.mps}}}
\]
Next, in every copy the interaction of the entangled state and the
measurement can be rewritten as shown.
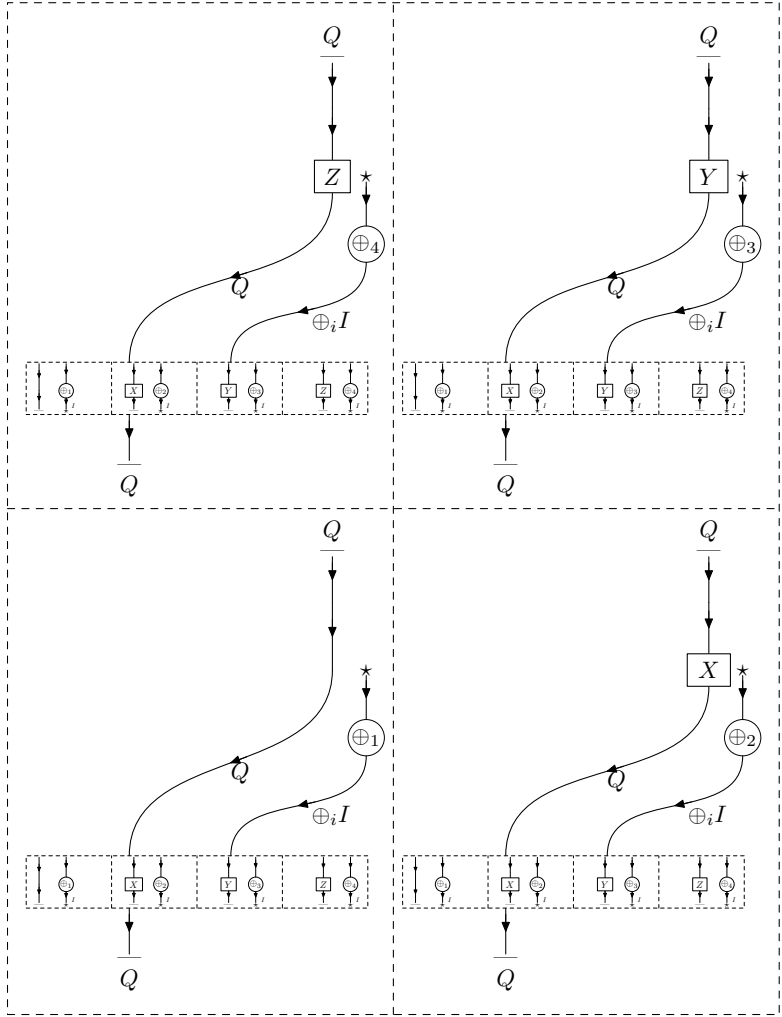
\[
 \raisebox{-0.5\height}{\scalebox{0.9}{
      \includegraphics{images/newproofnets_41.mps}}}
\]
Now we can open the box corresponding to Bob's non-determinism;  this
will leave us with sixteen copies of the system.  We won't draw all of
these copies, since  twelve of them can be erased: theses are the
cases the input that  
Bob is expecting does not match what Alice sends.  We are left with:
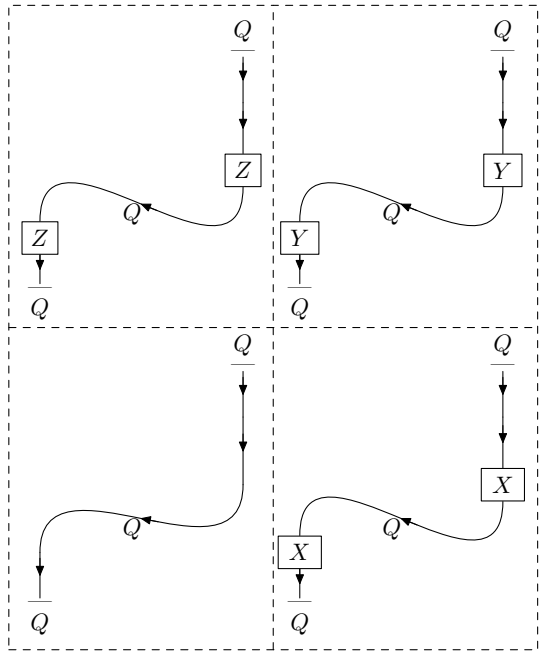
\[
\raisebox{-0.5\height}{\scalebox{0.9}{
      \includegraphics{images/newproofnets_42.mps}}}
\]
Now we simply note (and this is not a \emph{logical} axiom) that $X^2
= Y^2 = Z^2 = 1$ so we can simply remove these maps.  Hence we have
the normal form:
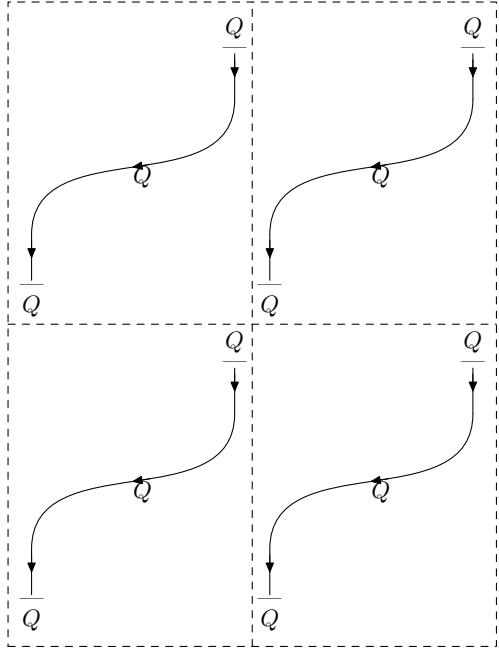
\[
\raisebox{-0.5\height}{\scalebox{0.9}{
      \includegraphics{images/newproofnets_43.mps}}}
\]
which show that in every possible world, Alice has successfully formed
a channel to Bob along which her state can be transmitted.

Although we presented the rewrites in the order that the steps of the
protocol would be carried out, in fact our proof-nets are strongly
normalising, so any order would produce the same results.

Having sketched the system, we move to the details.  First the
category theory, then in Sections~\ref{sec:tensor-sum-logic} and
\ref{sec:gener-proof-nets} the logic.

\section{Categorical Preliminaries}
\label{sec:categ-prel}

In this section, we introduce the necessary categorical structures,
compact closed  categories  and biproducts, 
and present their basic properties.  Much of this material is well known
so proofs are omitted.  Standard references are Mac Lane
\cite{MacLane:CatsWM:1971} and Kelly-Laplaza \cite{KelLap:comcl:1980};
other material is derived from the author's thesis \cite{Duncan:thesis:2006}.
Other sources are cited as needed.

Compact closed categories \cite{KelLap:comcl:1980} are abundant
throughout mathematics and computer science.  Examples include
\catRel, the category of sets and relations, finitely-generated
projective modules over a commutative ring and Conway games
(as categorified in \cite{Joyal:1977:conwaygames}).  
In $\mathbf{Hilb}$, the  category  of all Hilbert spaces, the
sub-category of  determined by the Hilbert-Schmitt maps is compact
closed, and more generally, the nuclear maps of any tensored *-category
\cite{ABP:nuclearideals:1999} for a compact closed sub-category.
Of course, \fdhilb, the category  of finite dimensional  Hilbert
spaces is compact closed.

In computer science compact closed  categories  have been studied in
the context of typed concurrency as interaction categories
\cite{AGN:1996:IntCats};  in logic, compact closed
categories are degenerate models of multiplicative linear logic 
\cite{AbrJag:NewFndGOI:1992,Loader:thesis:1994,HylSch:gomll:2003}; in
physics the category of $n$-dimensional  cobordisms, used in
topological quantum field, theory is compact closed
\cite{John-C.-Baez1995Higher-dimensio}.  More examples are easy to
find.

We define \emph{biproducts} and investigate their
basic properties in relation to compact closed categories.  Categories
with biproducts have been studied since the earliest days of
category theory as part of the theory of Abelian categories
\cite{Mitchell:ThCat:1965,MacLane:CatsWM:1971}.  Compact
closed categories with biproducts have been studied by Soloviev
\cite{Soloviev:superpositions:1987}.  A special case of compact closed
categories  with biproducts is \emph{Tannakien category}
\cite{Deligne:tannakiennes:1991}.  A recent contribution is by Houston
\cite{Houston2006Finite-Products}, who proved that every compact
closed category with products has biproducts.

\subsection{Monoidal Categories}
\label{sec:monoidal-categories}

\begin{definition}
  A category \catC is \emph{monoidal} if equipped with a functor
  $\otimes: \catC\times\catC\to\catC$, a distinguished neutral object
  $I$, and natural isomorphisms
  \begin{eqnarray*}
    \alpha_{A,B,C}:A\otimes(B\otimes C)\rTo^{\iso}(A\otimes B)\otimes C, \\
    \lambda_A:I\otimes A\rTo^{\iso} A, \qquad \rho_A:A\otimes I \rTo^{\iso} A.
  \end{eqnarray*}
  For the associativity morphism $\alpha$ we require that the pentagon 
  \begin{diagram}
    A\otimes(B\otimes(C\otimes D))&\rTo^\alpha&(A\otimes B)\otimes
    (C\otimes D)&\rTo^\alpha& ((A\otimes B)\otimes C)\otimes D \\
    \dTo^{\id{}\otimes\alpha} &&&& \uTo_{\alpha\otimes\id{}} \\
    A\otimes((B\otimes C)\otimes D)&&\rTo_\alpha&&(A\otimes(B\otimes C))\otimes D
  \end{diagram}
  commutes.  The isomorphisms $\lambda$ and $\rho$ express the
  neutrality of $I$; we require that the following diagram commutes:
  \begin{diagram}
    A\otimes(I\otimes B)&&\rTo^\alpha&&(A\otimes I)\otimes B\\
    &\rdTo^{\id{}\otimes\lambda}&&\ldTo_{\rho\otimes\id{}}&\\
    &&A\otimes B.&&
  \end{diagram}
\end{definition}

\begin{proposition}
  In a monoidal category the equality
  \[
  \lambda_I = \rho_I
  \]
  holds and the following diagrams commute:
  \begin{diagram}
    (A\otimes B)\otimes I & \rTo^\alpha & A\otimes(B\otimes I) 
    & \qquad &   
    (I\otimes A)\otimes B & \rTo^\alpha & I\otimes(A\otimes B) \\
    &  \rdTo(1,2)_\rho  \ldTo(1,2)_{\id{}\otimes\rho}& 
    & \qquad &   
    &  \rdTo(1,2)_{\lambda\otimes\id{}}  \ldTo(1,2)_\lambda & \\
    & A\otimes B & 
    & \qquad & 
    & A\otimes B. &
  \end{diagram}
  \begin{proof}See \cite{JS:1993:GeoTenCal1}.\end{proof}
\end{proposition}

\begin{definition}
  A monoidal category is \emph{symmetric}
 if it has a natural 
  isomorphism
\[
\sigma_{A,B}:A\otimes B\to B\otimes A
\]
such that
\begin{diagram}
  & (A\otimes B)\otimes C & \rTo^{\quad\sigma\quad} & C\otimes (A\otimes B) & \\
  \ldTo(1,2)^{\alpha^{-1}} &&&& \rdTo(1,2)^{\alpha}\\
  A\otimes(B\otimes C) &&&& (C\otimes A)\otimes B \\
  & \rdTo(1,2)_{\id{}\otimes\sigma} && \ldTo(1,2)_{\sigma\otimes\id{}}  & \\
  & A\otimes (C\otimes B) & \rTo_\alpha & (A\otimes C)\otimes B, & \\
\end{diagram}
and
\begin{diagram}
  A\otimes I &\rTo^\sigma & I\otimes A & \qquad & A\otimes B
  &\rTo^{\id{}} & A\otimes B \\
  & \rdTo(1,2)_\rho \ldTo(1,2)_\lambda & & \qquad & &
  \rdTo(1,2)_\sigma \ruTo(1,2)_\sigma & \\
  & A & & \qquad & & B\otimes A &
\end{diagram}
commute.
\end{definition}

\noindent
Mac Lane's celebrated coherence theorem states that any formal diagram
constructed from the $\alpha, \rho, \lambda$ and $\sigma$ will
commute.    A monoidal category is called \emph{strict}
if the isomorphisms $\alpha,\lambda$,and $\rho$ are all identities.
To minimise syntactic overhead we will make use of the
following theorem through this section:

\begin{theorem}[Mac Lane]\label{thm:maclaneStrict}
Every monoidal category $\catC$ is equivalent to some strict monoidal
category $\catA$.
\end{theorem}
\noindent  Note, however, that the category of proof-nets constructed in
Section~\ref{sec:gener-proof-nets} is \emph{not} strict: we will
produce non-trivial associativity and unit morphisms.

\subsection{Compact Closed Categories}
\label{sec:comp-clos-categ}
Let \catC be a symmetric monoidal category.  We say that \catC is \emph{compact
  closed} if every object $A$ has a chosen \emph{dual}%
\footnote{
  Some writers call this the ``adjoint'' in light of the relation
  between $A$ and $A^*$;  we use ``dual'' here to avoid confusion with
  the linear algebraic use of the word adjoint.}
 $A^*$ and maps 
\begin{eqnarray*}
  \eta_A: I\to A^*\otimes A \\
  \epsilon_A : A\otimes A^* \to I
\end{eqnarray*}
such that the composites
\begin{equation*}
  A \rTo^{\quad\id{A}\otimes\eta\quad} A\otimes A^* \otimes A
  \rTo^{\quad\epsilon \otimes \id{A}\quad} A
\end{equation*}
and
\begin{equation*}
A^* \rTo^{\quad\eta\otimes \id{A^*}\quad} A^*\otimes A\otimes A^* \rTo^{\quad\id{A^*}\otimes\epsilon\quad} A^*
\end{equation*}
are equal to $\id{A}$ and $\id{A^*}$ respectively.  We  call $\eta_A$
and $\epsilon_A$ the \emph{unit} and \emph{counit} maps.   

\begin{proposition}\label{prop:comclcanonicalisos}
In a compact closed  category we have natural isomorphisms:
\begin{eqnarray*}
  u : (A\otimes B)^* & \isomorphism & B^*\otimes A^* \\
  v : I^* & \isomorphism& I \\
  w : A^{**} & \isomorphism & A,
\end{eqnarray*}
\end{proposition}

\noindent A compact closed category which, in addition to being strictly
monoidal, has all of the isomorphisms $u,v,w$  equal to the identity
is called a \emph{strict compact closed category}.  Kelly and Laplaza
\cite{KelLap:comcl:1980} show that any compact closed category is
equivalent to a strict one, hence we will take the isomorphisms above
to be equalities whenever convenient.  

\begin{proposition}\label{prop:etaepsilonDNT}
  In a compact closed category the units and counits define dinatural
  transformations (see \cite{girard91normal})
  \begin{eqnarray*}
    && \eta:I \Rightarrow ((-)^*\otimes -) \\
    && \epsilon:(-\otimes  (-)^*) \Rightarrow I.
  \end{eqnarray*}
\end{proposition}

\noindent
We have a bijection between $\catC(A,B)$ and
$\catC(B^*,A^*)$:  given $f:A\to B$, define $f^*:B^*\to A^*$ by 
\begin{equation*}
\begin{diagram}
  B^* & \rTo^{\quad\eta_A \otimes \id{B^*}\quad} & A^*\otimes A\otimes B^*
  \\
  \dTo^{f^*} && \dTo_{\id{A^*}\otimes f \otimes\id{B^*}}
  \\
  A^* & \lTo_{\quad\id{A^*}\otimes\epsilon_{B}\quad} & A^*\otimes B\otimes B^*.
\end{diagram}  
\end{equation*}
We call $f^*$ the \emph{dual} of $f$.

\begin{proposition}\label{prop:star-is-functor}
  The operation $(-)^*$ defines a functor $\catC^{\text{op}}\to\catC$,
  which is an equivalence of categories.
  \begin{proof}
    We have $\id{A}^* = \id{A^*}$ immediately from the definition of
    dual, and $(f\circ g)^* = g^*\circ f^*$ follows from a routine
    calculation.  Taking \catC to be strict, we we have $A^{**} = A$,
    it follow from the defining property of compact closure that
    $f^{**} = f$, which gives the equivalence.
  \end{proof}
\end{proposition}

\noindent Since we have the equivalence between \catC and $\catC^{\text{op}}$,
any statement about some arrow applies equally well
to its dual.  In particular, results concerning units translate directly
into results about counits and vice-versa.  The duality of a compact
closed category gives a particularly strong 
form of monoidal closure.  Every arrow in the category has a
\emph{point} which represents it, and dually a \emph{copoint}.  These
representatives, the names and conames, will be crucial to our
treatment of entangled quantum states.

\begin{definition}\label{defn:nameconame}
  Let $f:A\to B$ in a compact closed category \catC.  Define the
  \emph{name} and \emph{coname} of $f$ to be the maps $\name{f}:I\to
  A^*\otimes B$ and $\coname{f}:A\otimes B^*\to
  I$ which are defined by the diagrams below.
  \begin{diagram}
    I & \rTo^{\eta_A} & A^*\otimes A & \qquad\qquad & A\otimes B^* &&
    \\
    & \rdTo_{\name{f}} &\dTo_{\id{A^*}\otimes f} &&
    \dTo{f\otimes\id{B^*}} &\rdTo^{\coname{f}} &\\
    && A^*\otimes B && B^*\otimes B &\rTo_{\epsilon_B} & I
  \end{diagram}
\end{definition}

\noindent An immediate consequence of this definition is the isomorphism of
hom-sets
\[
\catC(I,A^*\otimes B) \isomorphism \catC(A,B) \isomorphism
\catC(A\otimes B^*,I).
\]
\begin{lemma}\label{lem:absorption}
  Let \catC be compact closed and suppose we have arrows
  \[
  D\rTo^h A \rTo^f B \rTo^g C;
  \]
  then the following equations hold:
  \begin{align*}
    (\id{A^*}\otimes g ) \circ \name{f} &= \name{g\circ f}\;,\\
    ( h^*\otimes \id{B} ) \circ \name{f} &= \name{f\circ h}\;,\\    
\intertext{ and}
    (\coname{f}\otimes\id{C})\circ (\id{A}\otimes \name{g}) &= g\circ f
    \;.
  \end{align*}
\end{lemma}

\noindent We can also define partial versions of the name and coname; essentially currying and uncurrying.
\begin{lemma}[Partial Names and Conames]\label{lem:generalisednames}
In any compact closed category we have the following isomorphisms:
\begin{eqnarray}
  \catC(A\otimes C, B) & \isomorphism
  &
  \catC(A, C^*\otimes B) \label{eq:currying} \\
  \catC(A, C\otimes B) & \isomorphism
  &
  \catC( A \otimes C^*,B)
\end{eqnarray}
\begin{proof}
Since the two isomorphisms are dual, we prove only the first.
Define $F: \catC(A\otimes C, B)\to
\catC(A, C^*\otimes B)$
and $G :\catC(A, C^*\otimes B) \to \catC( A\otimes C, B)$ by 
\begin{eqnarray*}
  F:f \mapsto (\id{}\otimes f)\circ (\eta_A \otimes\id{}) \\
  G:g \mapsto (\epsilon_A\otimes\id{}) \circ (\id{}\otimes g) \\
\end{eqnarray*}
Their composition gives $GFf = (\epsilon_A\otimes\id{}) \circ (\id{}\otimes f) \circ (\id{}\otimes\eta_A\id{})$ from which
\begin{diagram}
  A\otimes C & 
  \rTo^{\id{}\otimes\eta_A \otimes\id{}} & 
  A \otimes A^* \otimes A \otimes C & 
  \rTo^{\id{}\otimes f} &
  A \otimes A^* \otimes B
\\
  & 
  \rdTo_{\id{}} &
  \dTo_{\epsilon_A\otimes\id{}} &
  &
  \dTo_{\epsilon_A\otimes\id{}}
\\
&&
  A\otimes C &
  \rTo_{f} &
  B
\end{diagram}
and hence $GF = \text{Id}$.  Similarly  $\text{Id} = FG$, which establishes the isomorphism.  
\end{proof}
\end{lemma}

\noindent Equation~(\ref{eq:currying}) essentially states that compact closed
categories are indeed closed with $B^A = A^*\otimes B$.  Since
$A^*\otimes I \isomorphism A^*$ this gives immediately the following.

\begin{corollary}\label{cor:comcl-are-starauto}
  Compact closed categories  are $*$-autonomous  \cite{Barr:*-AutoCats:1979}.
\end{corollary}

\noindent Hence compact closed  categories  are models of \MLL, and in
particular the linear $\lambda $-calculus;  albeit, these are rather strange
models equipped  with only one, self-dual, tensor.

\begin{definition}
  Let \catC be compact closed and define a map  
  \[
  \Tr : \catC( A \otimes C, B \otimes C) \to \catC (A,B)
  \]
  by setting 
  \[
  \Tr^C_{A,B}(f) = (\id{B} \otimes \epsilon_{C}) \circ 
           (f \otimes \id{C^*}) \circ
           (\id{A} \otimes \eta_{C^*})
  \]
  The map $\Tr(f)$ is called the \emph{trace of $f$}.
\end{definition}

\noindent
The trace so defined makes $\catC$ into a \emph{traced monoidal
  category} in the sense of Joyal, Street, and Verity
\cite{JSV:TraMonCat:1996}.  Few of the properties of the trace
are required here so we will not recapitulate the definition---in any case the relevant facts can be deduced from the
properties of $\eta$ and $\epsilon$.  We will, however, need the following
lemma:

\begin{lemma}[\cite{AHS:2002:GoILinComAlg}]
Suppose we have arrows $A\rTo^f B \rTo^g C$ in a symmetric traced
monoidal category; then: 
\[
g\circ f = \text{Tr}^B_{A,C}(\sigma_{B,C} \circ (f\otimes g))\;.
\]
\end{lemma}

\noindent The \emph{partial} trace defined above may be extended to a
full trace over any endomorphism $f:A\to A$ by setting
\[
\text{Tr}(f) = \text{Tr}_{I,I}^A( \rho\circ f \circ \rho^{-1}).
\]
In \fdhilb, the category of finite dimensional  Hilbert spaces, this
coincides with the usual trace, explicitly given by summing the
diagonal elements of a matrix representation of $f$.

\subsection{Scalars and Loops}
\label{sec:scalars-and-loops}

\begin{definition}
  In any monoidal category \catC the endomorphisms of the neutral
  element $\catC(I,I)$ are called the \emph{scalars}.
\end{definition}

\begin{lemma}\label{lem:scalarscommute}
The scalars form a commutative monoid with respect to composition.
\begin{proof}
Let $s,t\in\catC(I,I)$; then
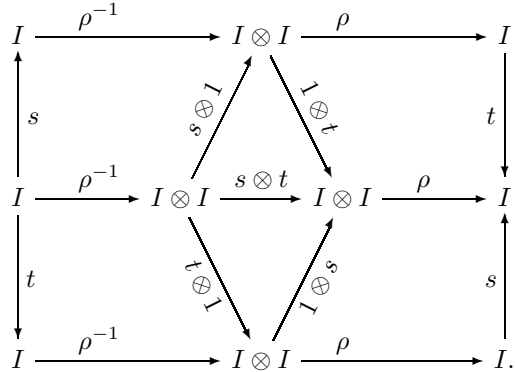
\begin{diagram}
  I & \rTo^{\rho^{-1}} && I\otimes I &\rTo^{\rho} && I \\
  \uTo_s && \ruTo(1,2)^{s\otimes\id{}} & & \rdTo(1,2)^{\id{}\otimes t} && \dTo^t \\
  I & \rTo^{\rho^{-1}} & I\otimes I & \rTo^{s\otimes t} & I\otimes I &  \rTo^{\rho} & I \\
  \dTo_t &&& \rdTo(1,2)_{t\otimes\id{}} \ruTo(1,2)_{\id{}\otimes s} &&& \uTo^s \\
  I & \rTo^{\rho^{-1}} && I\otimes I &\rTo^{\rho} && I.
\end{diagram}
\end{proof}
\end{lemma}

\begin{corollary}\label{cor:scalartensors}
  For scalars $s,t$ the composite
  \[
  I \isomorphism I\otimes I \rTo^{s\otimes t} I \otimes I \isomorphism  I
  \]
  is equal to $s\circ t = t\circ s$.
\end{corollary}

\begin{definition}\label{defn:scalarmultiple}
  Let \catC be a monoidal category.  Given a scalar $s$ and some arrow
  $f:A\to B$ define a scalar multiplication $s\bullet
  f$ by the
  composition:
  \begin{diagram}
    A &\rTo^{\rho^{-1}}& A\otimes I& \rTo^{f\otimes s} &B\otimes I &\rTo^\rho & B.
  \end{diagram}
\end{definition}

\noindent  
We could have defined $s\bullet f$ equivalently by multiplication on the 
  left rather than on the right as above.  Note that $u :=
  \lambda^{-1}\circ\rho$ is a natural isomorphism $(-\otimes
  I)\Rightarrow(I\otimes -)$, so the following diagram commutes
\begin{diagram}
    && A\otimes I& \rTo^{f\otimes s} &B\otimes I  & & \\
    A & \ruTo(2,1)^{\rho^{-1}}
    &\dTo^{u} & & \dTo^{u} &
    \rdTo(2,1)^\rho  & B\\
    &\rdTo(2,1)_{\lambda^{-1}} & I\otimes A& \rTo_{s\otimes f} &I\otimes B & \ruTo(2,1)_\lambda & 
\end{diagram}
and hence the two definitions coincide.

\begin{lemma}\label{lem:NTscalars}
  Each scalar $s$ determines a natural transformation
  $\text{Id}\Rightarrow \text{Id}$ such that $s\bullet f = f \circ s_A
  = s_B \circ f$.
  \begin{proof} The top and bottom edges define $s_A$ and $s_B$
    respectively:
    \begin{diagram}
      A & \rTo^{\rho^{-1}} & A\otimes I & \rTo^{\id{A}\otimes s} & A\otimes I
      & \rTo^\rho & A \\
      \dTo^f && \dTo^{f\otimes\id{I}} & \rdTo^{f\otimes s} &
      \dTo_{f\otimes\id{I}} && \dTo_f \\
      B & \rTo_{\rho^{-1}} & B\otimes I & \rTo_{\id{B}\otimes s} & B\otimes I
      & \rTo_\rho & B. \\
    \end{diagram}
    The outer squares commute due to naturality of $\rho$, and the
    middle due to the functoriality of the tensor.  Hence $s$ defines
    a natural transformation.  Note that the middle path from $A$ to
    $B$ is the definition of $s\bullet f$.
  \end{proof}
\end{lemma}

\begin{corollary}\label{cor:scalarequalities}
  The following are immediate.
  \begin{enumerate}
  \item $s\bullet (t\bullet f) = (s\circ t)\bullet f$
  \item $(s\bullet f)\circ(t\bullet g) = (s\circ t)\bullet(f\circ g)$
  \item $(s\bullet f)\otimes (t\bullet g) = (s\circ t)\bullet (f\otimes g)$
  \end{enumerate}
\end{corollary}

\begin{definition}\label{defn:dimension}
In a compact closed category \catC define the
\emph{dimension} of an object $A$, to be the following composite:
\[
\text{dim}_A = I \rTo^{\eta_A} A^*\otimes A \rTo^\sigma A\otimes A^*
\rTo^{\epsilon_A} I.
\]
\end{definition}

\noindent Of course, this is nothing more than the trace of
$\id{A}$.  The presence of these non-trivial scalars gives a
qualitative aspect even to freely constructed compact closed
categories.

\subsection{Freely Constructed Compact Closed Categories}
\label{sec:free-construction}

Define the set of \emph{endomorphisms} $E(\catA)$
by the disjoint union 
\[
E(\catA) = \sum_{A \in |\catA|} \catA(A,A),
\]
and let the set of \emph{loops} $[\catA]$ be the quotient of
$E(\catA)$ generated by the relation $f\circ g \sim g\circ f$ whenever
$A\rTo^f B \rTo^g A$. Let $\tau:E(\catA)\to [\catA]$ be the canonical
map onto the loops, and for each endomorphism $f$ write $[f]$ for
its image under $\tau$.  

The key theorem is the following of \cite{KelLap:comcl:1980}.
\begin{theorem}\label{thm:kellylaplazascalars}
Let 
\[
T : \catA^{\text{op}}\times \catA \times \catA^{\text{op}}\times \catA\times
\cdots \times \catA^{\text{op}}\times \catA \rTo \catB
\]
be a functor of $2n$ variables, let $K$ and $L$ be objects of \catB
and let $\alpha : K \Rightarrow T$ and $\beta:T\Rightarrow L$ be
natural transformations with typical components
\begin{gather}
  \alpha : K\rTo T(A_1,A_1,A_2,A_2,A_3,\ldots,A_{n-1},A_n,A_n)\;, \label{eq:KLloop-lemma-nt1}\\
  \beta :  T(B_1,B_2,B_2,B_3,\ldots,B_{n-1},B_n,B_n,B_1) \rTo L\;;\label{eq:KLloop-lemma-nt2}
\end{gather}
given maps 
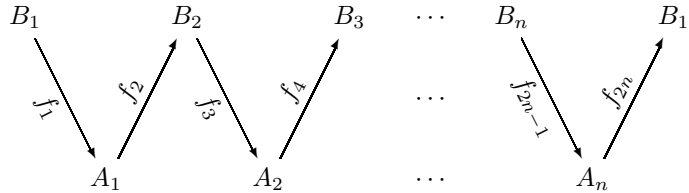
\begin{diagram}
  B_1 && B_2 && B_3 &\cdots&B_n && B_1 \\
&\rdTo(1,2)_{f_1}  \ruTo(1,2)^{f_2} &
&\rdTo(1,2)_{f_3}  \ruTo(1,2)^{f_4} &
& \cdots &
&\rdTo(1,2)_{f_{2n-1}}  \ruTo(1,2)^{f_{2n}} &
\\
& A_1 && A_2 && \cdots && A_n &&
\end{diagram}
the composite of (\ref{eq:KLloop-lemma-nt1}),
$T(f_1,f_2,f_3,\ldots,f_{2n-1},f_n)$ and (\ref{eq:KLloop-lemma-nt2}),
depends only on $[f_{2n}f_{2n-1}\cdots f_2f_1]$ so that $\alpha$
and $\beta$ give rise to a function $[\catA]\to \catB(K,L)$.
\end{theorem}

Taking $\catA = \catB$, $K = L = I$,  $\alpha = \eta$ and $\beta =
\epsilon$ gives a ready source of scalars in any compact closed
category \catA; indeed this is the dimension map given in Definition
\ref{defn:dimension}.  If the category is freely constructed these are the
\emph{only} non-trivial scalars.   This is a consequence of the more
general coherence theorem of Kelly and Laplaza.  Before stating the
theorem we must introduce some additional terminology, which will be
also be required later in this section.

\begin{definition}\label{defn:signed-set}
A \emph{signed set} $S$ is a function from a carrier set $|S|$ to the
set $\{+,-\}$. Given signed sets $R$ and $S$, let $R^*$ denote the
signed set with the opposite signing to $R$; let $R\otimes S$ be
the disjoint union of $R$ and $S$, such that $|R\otimes S| = |R|+|S|$.
\end{definition}

\begin{definition}
  An \emph{involution} is a category which is a coproduct of copies of
  the category  $\mathbf{2}$.  Given an involution $\sigma$, its
  object set $\sizeof{\sigma}$ can form a 
  signed set by assigning $-$ to the source and $+$ to the target of
  each arrow of $\mathbf{2}$.  Call $\sigma$ an \emph{involution on
    the signed set $S$} when this signing agrees with that of $S$.
\end{definition}

\noindent Given some category  \catA, we can construct the free compact closed
category  generated by \catA, which we call
$F\catA$. 
The objects of the $F\catA$ are constructed from those of $\catA$ by
repeated application of the  functors $-\otimes -$, $(-)^*$ and
the constant $I$. This characterisation may be used to inductively
construct a signed set $S(X)$ corresponding to each object $X$ of
$F\catA$.  Let 
\begin{align*}
  S(I) &= \emptyset\;, \\
  S(X\otimes Y) &= S(X) \otimes S(Y)\;,\\
  S(X^*) &= S(X)^*\;,\\
  S(A) &= \{ A \mapsto + \}\qquad\text{ if $A$ is an object of \catA}\;.
\end{align*}
The basic structure of arrows in $F\catA$ depends upon involutions on
the signed sets generated by its objects.

\begin{theorem}[Kelly-Laplaza]\label{thm:kellylaplaza}
Let \catA be a category;  each arrow $f:A\to B$ of the free compact
closed category generated by \catA is completely described by the following data:
\begin{enumerate}
\item An involution $\sigma$ on $S(A^*\otimes B)$;
\item A functor $\theta: \sigma \to \catA$ agreeing with $\sigma$ on
objects (i.e. a labelling of $\sigma$ with arrows of \catA);
\item A multiset $L$ of loops from \catA.
\end{enumerate}
\end{theorem}

\noindent The baroque statement of this theorem conceals its graphical
content.  One can view the objects of $F\catA$ as lists of positively
and negatively occurring objects of \catA, and an arrow between two
such lists is simply a collection of arcs, each connecting a negative
occurrence to a positive one, and labelled by an arrow of
$\catA$.
To compose arrows in $F\catA$ we simply connect
up the arcs, using the underlying composition in \catA.

From this point of view we can see an immediate limitation in the use
of such freely generated compact closed categories to model quantum
states.  Recall that when we interpret processes in categorical terms,
we view the objects as state spaces; hence the objects of the generating
category \catA are the state spaces of the elementary subsystems from
which our composite systems will be built.  A state of a compound
system, that is, an arrow $\psi : I \to X$ in $F\catA$, is thus
composed of \emph{pairs} of elementary systems related by some arrow
from $\catA$, and each pair is unconnected to the others.  
Quantum informatics attaches great importance to
\emph{entangled states}; that is, states which cannot be broken down
into their constituent parts.  However the above result states that
free compact closed categories can only result in bipartite
entanglement, which does not suffice to describe all entangled states.
To extend our reach we now introduce \emph{polycategories}.

\subsection{Compact Symmetric Polycategories}
\label{sec:comp-symm-polyc}

The reason that the free construction described above yields only
bipartite states is simple.  The states are based on the arrows of the
underlying category $\catA$, and an arrow has exactly two ends.  In
order to represent multipartite states we will need generators with
more than one and input and output, suggesting the need to construct
the compact closed  category  from a category which already has a
monoidal structure.  However the direct route leaves open the
problem of ensuring that the downstairs tensor (from the monoidal
category  \catA) and the upstairs tensors (freely generated in
$F\catA$) cohere correctly.  Worse, there is no reason to believe that
an arrow in a monoidal category is in any sense indecomposable among
its subsystems.

Fortunately there is a natural generalisation of category, a
\emph{polycategory}, whose arrows may have more than one object in
their domain and codomain.  The original notion of polycategory
\cite{Lambek1968Deductive-Syste,Szabo1975Polycategories} was 
introduced to study classical logic, where a sequent may have multiple
premises and conclusions; composition is defined by the cut-rule, so
one output is connected to one input.   Here we consider
\emph{compact symmetric} polycategories \cite{Duncan:thesis:2006},
where composition is defined by the \emph{multi-cut} rule, allowing
arbitrary vectors of inputs and outputs to be composed.  

\begin{definition}\label{def:polycat}
  A \emph{compact symmetric polycategory}, \catP, consists of a class 
  of objects \OBJC{P} and, to each pair $(\Gamma,\Delta)$ of finite
  sequences over \OBJC{P}, a set of polyarrows $\catP(\Gamma,\Delta)$.
  Given a non-empty sequence of objects $\Theta$ and 
  poly-arrows
  \[
  \Gamma\rTo^f \Delta_1,\Theta, \Delta_2 \qquad \text{ and } \qquad
  \Gamma_1,\Theta,\Gamma_2\rTo^g \Delta
  \]
  we may form the composition%
  \index{$\polycomp{i}{j}{k}$}
  \[
  \Gamma_1,\Gamma,\Gamma_2 \rTo^{g \polycomp{i}{j}{k} f} \Delta_1,\Delta,\Delta_2
  \]
  where $\sizeof{\Delta_1} = i$, $\sizeof{\Gamma_1} = j$ and
  $\sizeof{\Theta} = k > 0$. For each object $A$ there is an identity
  arrow $\id{A}: \langle A \rangle \to \langle A \rangle$ for the
  singleton sequence $\langle A \rangle$.
\end{definition}

In general there are many ways to compose the  polyarrows, and many
equations which must be satisfied.  We will spare the reader the full
definition\footnote{For the full glory of its coherence equations,
  and also the proof of the  theorem cited, see \cite{Duncan:thesis:2006}.},
and instead we will offer a theorem in the spirit of the Kelly-Laplaza
result cited above, characterising the free compact closed category
generated by such a polycategory.  Before proceding we note the
most important point about these polycategories:  there is no nullary
composition and no tensor product.  Each input of a polyarrow
has a path (not necessarily directed) to each output, and hence
despite having many inputs and outputs, polyarrows cannot be decomposed
into non-interacting parts.

Before we can state the representation theorem we must make some
definitions.

\begin{definition}\label{defn:graph}
  A \emph{graph} consists of a 5-tuple $(V,E,C,s,t)$ where $V$,$E$,
  and $C$ are sets, respectively of \emph{vertices}, \emph{edges}, and
  \emph{circles}, and $s$ and $t$ are maps
  \begin{diagram}
    E & \pile{ \rTo^s \\ \rTo_t} & V
  \end{diagram}
  which we call \emph{source} and \emph{target}.  Let  $\text{in}(v)$ and
  $\text{out}(v)$ be $V$-indexed subsets of $E$ defined by
  \begin{align*}
    \text{in}(v) & = t^{-1}(v)\\
    \text{out}(v) & = s^{-1}(v).
  \end{align*}
  The \emph{in-degree} of a vertex $v$ is the
  cardinality of $\text{in}(v)$%
\index{$\text{in}(\cdot)$} and the \emph{out-degree} is the cardinality
  of $\text{out}(v)$%
\index{$\text{out}(\cdot)$}.  The
  \emph{degree} of a vertex is the sum of its in-
  and out-degrees.
\end{definition}

\begin{definition}\label{defn:opengraph}
A \emph{open graph} is a pair $(G,\partial G)$ of an
underlying graph $G = (V,E,C,s,t)$ and a
distinguished subset of the degree one vertices $\partial
G$%
\index{$\partial G$, boundary of graph $G$} called
the \emph{boundary} of $G$; $V-{\partial G}$ is called the
\emph{interior} of $G$, written $I_G$.  If a vertex $x \in \partial G$
it is an \emph{outer} or \emph{boundary node}; otherwise it is an
\emph{inner} or \emph{interior node}. 
\end{definition}

\begin{definition}\label{defn:circuit}
  A \emph{circuit} $\Gamma = (G, \dom\Gamma, \cod\Gamma,
  <_{\text{in}(\cdot)}, <_{\text{out}(\cdot)})$ where:
  \begin{itemize}
  \item $G = ((V,E,C,s,t),\partial G)$ is an open graph;
  \item $\dom\Gamma$%
\index{$\dom\Gamma$} and $\cod\Gamma$%
\index{$\cod\Gamma$} are totally ordered sets such
    that  $\partial G =  \dom\Gamma + \cod\Gamma$;
  \item $<_{\text{in}(\cdot)}$%
\index{$<_{\text{in}(\cdot)}$} is a family of maps, indexed by $V$ such
    that 
    \[
    <_{\text{in}(v)} : \text{in}(v) \rTo^\isomorphism \mathbb{N}_k
    \]
    where $k = \sizeof{\text{in}(v)}$.
  \item $<_{\text{out}(\cdot)}$%
\index{$<_{\text{out}(\cdot)}$} is a family of maps, indexed by $V$ such
    that 
    \[
    <_{\text{out}(v)} : \text{out}(v) \rTo^\isomorphism \mathbb{N}_{k'}
    \]
    where $k' = \sizeof{\text{out}(v)}$.
\end{itemize}
\end{definition}

\begin{figure}[htbp]
  \centering
  \includegraphics{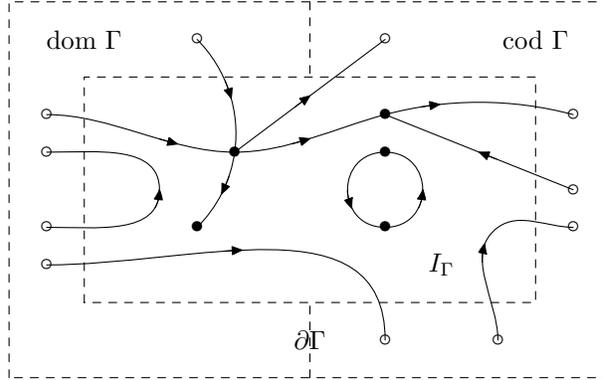}
  \caption{Anatomy of a circuit}
  \label{fig:circuitanatomy}
\end{figure}

\noindent As suggested by their name, the purpose of the two maps
$<_{\text{in}(\cdot)}$ and $<_{\text{out}(\cdot)}$ is to impose a
linear order on $\text{in}(v)$ and $\text{out}(v)$.  Since the maps
give a bijective correspondence between $\text{in}(v),\text{out}(v)$
and an initial segment of the naturals, the order in $\mathbb{N}$
lifts, and hence we will often simply treat these sets are ordered,
and write $<$ for this ordering whenever unambiguous to do so.

For simplicity, in the following we will consider a polycategory \catP
which is freely generated from some set of basic arrows\footnote{This
  is not essential, but will greatly simplify the subsequent
  discussion of generalised proof-nets;  see the discussion of
  homotopy in \cite{Duncan:thesis:2006} for the details.} called \ARRC{P}.

\begin{definition}\label{def:labelling}
Given a polycategory \catP , an \emph{\catP-labelling} for a circuit
$\Gamma$ is a pair of maps $\theta = (\theta_O,\theta_A)$ where
\begin{gather*}
  \theta_O : E\!+\!C\rTo \OBJC{P}  \\
  \theta_A : V \rTo \ARRC{P} 
\end{gather*}
such that for each node $f$, $\text{in}(f) = \langle a_1, \ldots, a_n\rangle$
and $\text{out}(f) = \langle b_1, \ldots, b_m\rangle$ imply 
\begin{eqnarray*}
\dom(\theta f) = \theta a_1, \ldots, \theta a_n \\
\cod (\theta f) = \theta b_1, \ldots, \theta b_m,
\end{eqnarray*}
and subject to the further restriction that $\theta_A(v) = \id{A}$ if
and only if $v\in {\partial \Gamma}$.
Call a circuit $\Gamma$ \emph{\catP-labellable} if there exists an 
\catP-labelling for it; if $\theta$ is a labelling for $\Gamma$, then
 the pair $(\Gamma,\theta)$ is an \emph{\catP-labelled circuit}.
\end{definition}

The boundary nodes perform a different role to the interior nodes.  
The  incidence of the unique edge at a boundary vertex $b$
defines a signing on the boundary: we say that $b$ is positive if it
has in-degree 1; and negative if its out-degree is 1.  We
will take the boundary vertices as labelled by objects of \catP rather
than the corresponding identity maps, and treat the boundary as an
\OBJ{P}-labelled signed set.

We denote the class of \catP-labelled circuits \CirP; it forms a
monoidal category in a rather natural way.  The objects of \CirP are
signed sequences of objects from \OBJ{P}.  An arrow
from $f: A \to B$ is defined by a \catP-labelled circuit whose
codomain is $B$ and whose domain is $A^*$ (i.e $A$ with the opposite
signing).  Composition is defined by joining two circuits at their
respective domain and codomain vertices, and erasing the vertices.  The tensor
product can be defined by taking the disjoint union of the circuits,
and concatenating the domain and codomains.   We leave the reader to
fill in the details.  We can now state the promised representation
theorem:

\begin{theorem}
  \CirP is the free compact closed category generated by the compact
  symmetric polycategory \catP.
\end{theorem}

\noindent 
Since categories are a special case of polycategory (where all the
arrows are between singleton sequences) we can ask: what is \CirP when
\catP is just a normal category?  In this case, \CirP is exactly
the same as $F\catP$ as per Kelly-Laplaza.

One way to understand the generalisation in going from a category of
generators to a polycategory of generators is by considering the case with
only one ground type.  If we have a category  with just one
object, we can view its arrows as evolutions of this state space.  On
the other hand, if we have a polycategory with a single object, the arrows
are in some sense \emph{interactions} between systems of that type,
possibly fusing or splitting, producing a different number of systems
than began the interaction.  By moving to the compact closed  category
generated by a polycategory of interactions we avoid the restriction
to bipartite states mentioned earlier.

This concludes the multiplicative structures, now we move onto the
additives.

\subsection{Zero Objects}
\label{sec:zero-objects}

\begin{definition}\label{defn:zero-object}
  In any category  \catC a \emph{zero object} is an object, denoted
  \ZERO, which is both initial and terminal.
\end{definition}

\noindent By its initiality, there is a unique map from \ZERO to every object,
and dually there is a unique map from each object to \ZERO.  Hence there is
unique map 
\begin{diagram}
A & \rTo^{\qquad} & \ZERO & \rTo^{\qquad} &B   
\end{diagram}
between every pair of objects $A$ and $B$.  This map is called the
\emph{zero map} and denoted $0_{A,B}$.  Since \ZERO is initial and
terminal any map composed with a zero map is again a zero map; the
zeros form a two-sided ideal with respect to composition among the
arrows \catC.  Hence the following diagram commutes:
\begin{diagram}
  A &\rTo & \ZERO & \rTo & B \\
  \dTo<f & & \dTo>{\id{\ZERO}} && \dTo>g\\
  C &\rTo & \ZERO & \rTo & D,
\end{diagram}
which makes the family $0_{A,B}$  natural in both $A$ and $B$.

A useful family of arrows in a category  with \ZERO is the Kronecker
delta $\delta_{ij}:A_i\rTo A_j$, defined for all pairs of objects
$A_i,A_j$ as
\[
\delta_{ii} = \id{A_i} \qquad \delta_{ij} = 0_{A_i,A_j}
\]

\begin{lemma}
  If $\id{A} = 0_{A,A}$ then $A \isomorphism  \ZERO$.
\end{lemma}
\begin{proof}
  Note that the composite
  \begin{diagram}
    \ZERO & \rTo^{!_A} & A & \rTo^{!^A} & \ZERO  
  \end{diagram}
  is equal to $0_{\ZERO,\ZERO}$, which by uniqueness is equal to
  $\id{\ZERO}$.  Thus ${!_A}\circ{!^A} = \id{A}$ and ${!^A}\circ{!_A}
  = \id{\ZERO}$, which gives the isomorphism.
\end{proof}

\begin{proposition}\label{prop:zero-tensor-A-is-zero}
  Let \catC be a monoidal closed category with a zero
  object.  Then $A\otimes\ZERO \isomorphism \ZERO$.
\end{proposition}
\begin{proof}
  Since \catC is closed,
  \[
  \catC(A\otimes\ZERO,B) \,\iso \,\catC(\ZERO, A\multimap B) \,\iso \, \{ * \}.
  \]
  Taking $B = A \otimes \ZERO$ implies  $\id{A\otimes\ZERO} = 0$, and
  hence the result follows.
\end{proof}

\noindent Monoidal closure is required;  if we take the tensor to be a
coproduct, e.g direct sum of vector spaces, it is clear that the
isomorphism does not hold. 

\begin{corollary}
  Given $f:A\to B$ an arrow of \catC, $f\otimes 0_{C,D} = 0_{A\otimes
    C, B\otimes D}$. 
\end{corollary}

\noindent
If the zero object is also the neutral object for the tensor, then the
entire category collapses to a single object via $A \isomorphism
A\otimes \ZERO \isomorphism  \ZERO$.  So any Cartesian closed
category with zero is trivial.  Note that in a compact closed category
 \catC 
with a terminal object $\mathbf{1}$, by duality  $\mathbf{1}^*$ is initial.
If the terminal object is the monoidal unit then the isomorphism $I
\isomorphism I^*$ makes $I$ the zero object, and hence the category
collapses.

\begin{proposition}
  If \catC is compact closed  with respect to a product then it is trivial.
\end{proposition}

\subsection{Biproducts}
\label{sec:biproducts}

In any category \catC with finite products and coproducts every map 
\[
\coprod A_i \rTo^f \prod A_j
\]
has a ``matrix'' representation $(f_{ij})$ where each $f_{ij}$ is given by
the composite 
\begin{diagram}
  A_i & \rTo^{\mathrm{in}_i} & \coprod A_i & \rTo^f & \prod A_j &
  \rTo^{\pi_j} & A_j
\end{diagram}
with $\mathrm{in}_i$ and $\pi_j$ the appropriate injections and
projections.  Supposing that \catC also has a zero object there
is canonical map $\mathbb{1}:\coprod A_i \rTo \prod A_i$ whose
matrix is the identity $\mathbb{1} = (\delta_{ij})$.

\begin{definition}\label{defn:biproduct}
  A category \catC has finite biproducts if it has finite products and
  coproducts, such that 
  \begin{itemize}
  \item the unique map $0\rTo 1$ is invertible; and
  \item the canonical map $\mathbb{1} : A\coprod B \rTo A\prod B$ is an
  isomorphism  for all objects $A, B$.
  \end{itemize}
\end{definition}
If \catC has biproducts, for all objects $A$ and $B$, there is a
unique (upto isomorphism) object $A\oplus B$ and maps 
\begin{equation}
\begin{diagram}
  A & \pile{\rTo^{\textrm{in}_1} \\ \lTo_{\pi_1}} & 
  A\oplus B & \pile{ \lTo^{\textrm{in}_2} \\ \rTo_{\pi_2}} & B
\end{diagram}\label{eq:biproduct-diagram}
\end{equation}
such that  $(A\oplus B, \pi_1,\pi_2)$ is a product and $(A\oplus
B,\mathrm{in}_1,\mathrm{in}_2)$ is a coproduct.  A choice of $A\oplus
B$ for every pair of objects makes $-\oplus -$ into a functor
$\catC\times \catC \rTo \catC$ whose action on arrows $f_1\oplus f_2$
is given by
\[
\pi_i \circ (f_1\oplus f_2) = f_i \circ \pi_i 
\qquad\text{ for } i =
1,2
\]
or alternatively
\[
(f_1\oplus f_2) \circ \mathrm{in}_i = \mathrm{in}_i \circ  f_i 
\qquad\text{ for } i = 1,2.
\]

\begin{lemma}
  In a category with biproducts we have the following natural
  isomorphisms:
  \begin{itemize}
  \item $(A\oplus B)\oplus C \isomorphism A \oplus (B\oplus C)$ ;
  \item $A\oplus B \isomorphism  B\oplus A$ ;
  \item $A\oplus\ZERO \isomorphism A  \isomorphism  \ZERO\oplus A$.
  \end{itemize}
\end{lemma}
\begin{proof}
  All these isomorphisms hold for products (and also coproducts) hence
  they hold for the biproduct.
\end{proof}

\noindent
We have natural diagonal and codiagonal maps,
\begin{gather*}
  \Delta_A : A \rTo A\oplus A\\
  \nabla_A: A\oplus A \rTo A 
\end{gather*}
defined as
\[
\Delta_A = \langle \id{A},\id{A}\rangle \qquad \nabla_A = [\id{A},\id{A}].
\]
It is useful to note the equations
\begin{gather*}
  (f\oplus g)\circ \Delta_A = \langle f,g\rangle, \\ 
  \nabla_A \circ (f\oplus g)= [f,g].
\end{gather*}

\begin{definition}
\label{defn:addition-of-arrows}  Let $f,g:A \rTo B$; then define $f+g$ as the composite
  \begin{diagram}
    A & \rTo^{\Delta_A} & A\oplus A& \rTo^{f\oplus g} & B\oplus B &
    \rTo^{\nabla_B}& B.
  \end{diagram}
\end{definition}

\begin{proposition}\label{prop:biproduct-implies-semi-additive}
  In a category with biproducts \catC, the addition of
  Definition~\ref{defn:addition-of-arrows}:
  \begin{itemize}
  \item makes each hom-set $\catC(A,B)$ is a commutative monoid; and 
  \item distributes over composition.
  \end{itemize}
\end{proposition}
\begin{proof}
  The addition is associative due to the following diagram
\[
\scalebox{0.725}{
  \begin{diagram}
    && A\oplus A & \rTo^{\Delta_A\oplus \id{A}} & (A\oplus A)\oplus A
    & \rTo^{(f\oplus g)\oplus h} & (B\oplus B)\oplus B & \rTo^{\nabla_B
    \oplus \id{B}} &B\oplus B && \\
    A & \ruTo(2,1)^{\Delta_A}&&& \dTo<\isomorphism  && \dTo>\isomorphism  &&&\rdTo(2,1)^{\nabla_B}&B \\
    &\rdTo(2,1)_{\Delta_A}& A\oplus A & \rTo_{\id{A}\oplus\Delta_A} & A\oplus (A\oplus A)
    & \rTo_{f\oplus (g\oplus h)} & B\oplus (B\oplus B) & \rTo_{\id{B}\oplus\nabla_B} &B\oplus B &\ruTo(2,1)_{\nabla_B} &
  \end{diagram}
}
\]
  and commutative since 
  \begin{diagram}
    && A\oplus A & \rTo^{f\oplus g} & B\oplus B &&\\
    A & \ruTo(2,1)^{\Delta_A} & \dTo<\sigma && \dTo>\sigma &
    \rdTo(2,1)^{\nabla_B} & B \\
    &\rdTo(2,1)_{\Delta_A} & A\oplus A & \rTo_{g\oplus f} & B\oplus B
    & \ruTo(2,1)_{\nabla_B} &
  \end{diagram}
  commutes.  The neutral element of $\catC(A,B)$ is $0_{A,B}$ since
  the following commutes:
  \begin{diagram}
    & A && \rTo^f && B & \\
    \ldTo(1,2)<{\Delta_A} &&\rdTo(1,2)>{\isomorphism} &&
    \ruTo(1,2)<{\isomorphism  } &&\luTo(1,2)>{\nabla_B} \\
    A \oplus A & \rTo_{\id{A}\oplus 0} & A\oplus \ZERO & \rTo_{f\oplus
        0}& B\oplus \ZERO &\rTo_{\id{B}\oplus 0}& B\oplus B
  \end{diagram}
  hence \catC is enriched over commutative monoids.  To see that the
  addition distributes over composition recall the identities 
  \begin{gather*}
    \langle f,g\rangle \circ h = \langle f\circ h, g\circ h\rangle \\
    h \circ [f,g] = [h\circ f, h\circ g]
  \end{gather*}
  hence 
  \begin{align*}    
    g\circ (f + f') \circ h &= g\circ \nabla \circ(f \oplus  f')
    \Delta \circ h \\
    &= [g,g]\circ (f\oplus f') \circ \langle h,h\rangle \\
    &= \nabla \circ ((g\circ f\circ h) \oplus(g\circ f'\circ h) )
    \circ \Delta \\
    & = (g\circ f\circ h) + (g\circ f'\circ h).
  \end{align*}
\end{proof}

\begin{proposition}\label{prop:biproduct-injections-projecion}
In a category with biproducts the injections and projections shown in
Eq~(\ref{eq:biproduct-diagram}) satisfy
\begin{gather*}
  \pi_i \circ \mathrm{in}_j = \delta_{ij} \qquad \text{ for } i,j= 1,2\\
   \mathrm{in}_1\circ \pi_1 + \mathrm{in}_2\circ \pi_2 = \id{A_1\oplus
     A_2}.
\end{gather*}
\end{proposition}
\begin{proof}
  For any product $A\times B$ we have $\pi_1\times \pi_2 \circ \Delta = \langle
  \pi_1, \pi_2 \rangle = \id{A\times B}$  and dually for any coproduct
  $\nabla\circ\mathrm{in}_1+\mathrm{in}_2 = \id{A+B}$.  Hence
  \[
  \mathrm{in}_1\circ \pi_1 + \mathrm{in}_2\circ \pi_2 =
  \nabla \circ  (\mathrm{in}_1 \oplus \mathrm{in}_2) \circ (\pi_1
  \oplus \pi_2) \circ \Delta = \id{A\oplus B}.
  \]
  Due to the universal property of the biproduct, the canonical map
  from $A\oplus B$ to itself is equal to $\id{A\oplus B}$.  Therefore 
  \[
  \pi_i \circ \mathrm{in}_j = \pi_i \circ \mathbb{1}\circ
  \mathrm{in}_j = \delta_{ij}.
  \]
\end{proof}

The binary biproduct may be generalised to arbitrary finite families
of objects $A_1,\ldots,A_n$ by
iteration.  Upto an associativity isomorphism, the $n$-fold biproduct
is characterised by the diagram
\begin{diagram}
  A_i & \rTo^{\mathrm{in_i}} &\bigoplus_k A_k & \rTo^{\pi_j} & A_j 
\end{diagram}
subject to the equations
\begin{gather*}
  \pi_i\circ \mathrm{in}_j = \delta_{ij},\\
  \sum_k \pi_k \circ \emph{in}_k = \id{\oplus_k A_k}.
\end{gather*}
Arrows between biproducts have matrix representations as described at
the start of this section, and composition of of arrows gives the
usual matrix multiplication.
\begin{proposition}\label{prop:matrix-calculus-for-biproducts}
  Suppose we have arrows
  \begin{diagram}
    \bigoplus_i A_i & \pile{\rTo^f\\\rTo_{f'}} & \bigoplus_j B_j & \rTo^g & \bigoplus_k C_k
  \end{diagram}
  then $(g\circ f)_{ik} = \sum_j (g_{jk}\circ f_{ij})$ and
  $(f+f')_{ij} = f_{ij} +f'_{ij}$.
\end{proposition}
\begin{proof}
  Let $h = g\circ f$ ; then
  \begin{align*}
   h_{ik} & = \pi_k\circ g \circ f \circ \mathrm{in}_i \\ &=  \pi_k\circ g
   \circ \id{\oplus B_j}\circ f \circ \mathrm{in}_i \\
   &=  \pi_k\circ g \circ (\sum_j \pi_j \circ \mathrm{in}_j) \circ f
   \circ \mathrm{in}_i\\
   &= \sum_j \pi_k\circ g \circ \mathrm{in}_j \circ \pi_j \circ f \circ
   \mathrm{in}_i \\
   &= \sum_j g_{jk}\circ f_{ij}.
  \end{align*}
  The second equation follows directly from the naturality of the
  diagonal and codiagonal maps.
\end{proof}

\noindent
Proposition \ref{prop:biproduct-implies-semi-additive} implies that
any category with biproducts is enriched over $\mathbf{CMon}$, the
category  of commutative monoids;  conversely, we have the following.

\begin{proposition}\label{prop:CMon-inplies-bips}
  Let \catC be a \textbf{CMon}-category with a \ZERO object and, for
  every pair of objects $A$ and $B$,  a diagram
  (\ref{eq:biproduct-diagram})  such that proposition
  \ref{prop:biproduct-injections-projecion} holds;  then \catC has
  biproducts.
\end{proposition}

\begin{theorem}
  A $\mathbf{CMon}$-category has products (or coproducts) if and only
  if it has biproducts.
\end{theorem}
\begin{proof}
  This is a fairly trivial modification of Mac Lane
  \cite{MacLane:CatsWM:1971} VIII.2, Theorem~2.
\end{proof}

\begin{definition}\label{defn:semi-additive-functor}  Call a
  $\mathbf{CMon}$-category  \emph{semi-additive} if it has \ZERO and a
  biproduct for each pair of its objects.
  Let \catA and \catB be  $\mathbf{CMon}$-categories with zero objects; a
  functor $F:\catA \rTo \catB$ is \emph{semi-additive} if $F\ZERO =
  \ZERO$ and $Ff + Fg = F(f+g)$ for all parallel arrows $f,g$ in \catA.
\end{definition}
\begin{proposition}
  Let \catA have biproducts and let \catB be $\mathbf{CMon}$-enriched
  with \ZERO;  then a functor $F\catA \to \catB$ is semi-additive if
  and only if it carries every biproduct diagram in \catA to a
  biproduct diagram in \catB.
\end{proposition}
\begin{proof}
  See Mac Lane \cite{MacLane:CatsWM:1971} VIII.2, Proposition 4.
\end{proof}

\noindent
Given a category \catC we can construct the free biproduct
structure on \catC, first by freely enriching \catC over
$\mathbf{CMon}$ and then taking matrices over the resulting
category.

\begin{proposition}\label{prop:freeCMon}
  Let $\catC_{\mathbb{N}}$ be the category whose objects are those of
  \catC and where $\catC_{\mathbb{N}}(A,B) = \mathbb{N}(\catC(A,B))$,
  the free commutative monoid on $\catC(A,B)$.  Then $\catC_{\mathbb
    N}$ is $\mathbf{CMon}$ enriched and the inclusion of
  \catC into $\catC_{\mathbb N}$ is a universal arrow from \catC to a
  $\mathbf{CMon}$-category.
\end{proposition}

\begin{proposition}\label{prop:freeMatrix}
  Let \catC be a   $\mathbf{CMon}$-category and let $\mathbf{Matr}(
  \catC)$ be the category whose objects are $n$-tuples of objects of
  \catC, for $n \geq 1$, and whose arrows are matrices of arrows
  \catC.  Then $\mathbf{Matr}(\catC)$  is semi-additive, and the
  evident semi-additive embedding of \catC into $\mathbf{Matr}(
  \catC)$ is universal among semi-additive functors from \catC to
  semi-additive categories. 
\end{proposition}

\subsection{Compact Closed Categories with Biproducts}
\label{sec:strongly-comp-clos-w-bips}

\begin{proposition}\label{prop:otimes-distributes-over-oplus}
  Let \catC be a monoidal closed category  with biproducts;  then there are
  natural distribution isomorphisms
  \begin{gather*}
    A\otimes (B\oplus C) \isomorphism  (A\otimes B)\oplus (A\otimes
    C)\\
    (A\oplus  B)\otimes C \isomorphism  (A\otimes C)\oplus (B\otimes
    C)
  \end{gather*}
\end{proposition}
\begin{proof}
  Since $A\otimes -$ is a left adjoint it preserves colimits and
  hence the diagram 
  \begin{diagram}
    A\otimes B & \rTo^{\quad\id{A}\otimes \mathrm{in}_1\quad} &
    A\otimes (B\oplus C) &
    \lTo^{\quad\id{A}\otimes \mathrm{in}_2\quad} & A\otimes C 
  \end{diagram}
  is a coproduct and hence $A\otimes (B\oplus C) \isomorphism
  (A\otimes B)\oplus (A\otimes C)$.  The right hand distribution is
  similar.
\end{proof}
\begin{corollary}
  In a monoidal closed category \catC with biproducts, the functor $A\otimes
  - : \catC\to\catC$ is additive.
\end{corollary}

\noindent 
In fact we can easily construct the distribution isomorphisms
explicitly.  Let 
\[
d_{A,B,C} = \langle \id{A}\otimes \pi_1, \id{A}\otimes \pi_2\rangle
\]
and 
\[
d^{-1}_{A,B,C} = [ \id{A}\otimes \mathrm{in}_1, \id{A}\otimes \mathrm{in}_2].
\]
Then 
\begin{multline*}
    [ \id{A}\otimes \mathrm{in}_1, \id{A}\otimes \mathrm{in}_2]\circ
    \langle \id{A}\otimes \pi_1, \id{A}\otimes \pi_2\rangle 
    \\
      \begin{split}
    \qquad & = 
    \nabla \circ 
     ((\id{A}\otimes \mathrm{in}_1) \oplus(\id{A}\otimes \mathrm{in}_2))
    \circ 
    ((\id{A}\otimes \pi_1)\oplus (\id{A}\otimes \pi_C)) \circ 
    \Delta
    \\
    \qquad & = \nabla\circ(\id{A}\otimes (\mathrm{in}_1\circ\pi_1)) \oplus
    (\id{A}\otimes (\mathrm{in}_2\circ\pi_2))\circ \Delta
    \\
    \qquad & = (\id{A}\otimes (\mathrm{in}_1\circ\pi_1)) + (\id{A}\otimes
    (\mathrm{in}_2\circ\pi_2))
    \\
    \qquad & = \id{A}\otimes(\mathrm{in}_1\circ\pi_1 +
    \mathrm{in}_2\circ\pi_2)
    \\
    \qquad & = \id{A}\otimes \id{B\oplus C} = \id{A\otimes (B\oplus C)}
  \end{split}
\end{multline*}

\noindent If a compact closed category \catC has a binary product $ - \times -$ then the
duality $(\cdot)^*$ sends every product diagram
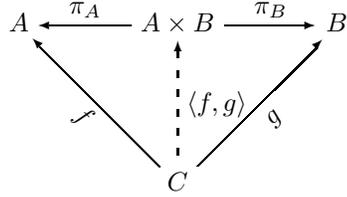
\begin{diagram}
  A & \lTo^{\pi_A\quad} & A\times B & \rTo^{\quad\pi_B} & B \\
  & \luTo<f & \uDashto>{\langle f, g \rangle} & \ruTo>g \\
  && C &&
\end{diagram}
to a coproduct diagram
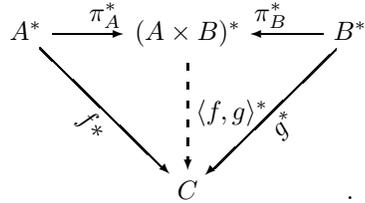
\begin{diagram}
  A^* & \rTo^{\pi^*_A} & (A\times B)^* & \lTo^{\pi^*_B} & B^* \\
  & \rdTo<{f*} & \dDashto>{\langle f, g \rangle^*} & \ldTo>{g^*} \\
  && C &&.
\end{diagram}
As mentioned earlier, if \catC has a terminal object $\mathbf{1}$ then
$\mathbf{1}^*$ is initial.  Hence the question of whether or not \catC
has biproducts boils down to whether the canonical maps
\begin{gather*}
\ZERO \rTo \mathbf{1} \\
(A\times B)^* \rTo A^* \times B^*  
\end{gather*}
are isomorphisms.  It turns out that this is always the case.

\begin{proposition}{(Houston)}
  If a compact closed  category  \catC has all finite products (or
  coproducts) it has all finite biproducts.
\end{proposition}
\begin{proof}
  See Houston \cite{Houston2006Finite-Products}.
\end{proof}

\begin{corollary}\label{cor:star-is-bip-nice}
  In any compact closed category with biproducts:
  \begin{itemize}
  \item  we have natural isomorphisms
    \[
    \ZERO \isomorphism  \ZERO^* \qquad\qquad (A\oplus B)^*\isomorphism
    A^*\oplus B^* ; 
    \]
  \item the duality $(\cdot)^*$ is an additive functor.
  \end{itemize}
\end{corollary}

\noindent 
It then follows that we may choose the biproduct in any compact closed
category so that the equation
\[
\mathrm{in_A}^* = \pi_{A^*}
\]
holds for all objects $A$. 

We now turn our attention to the construction of the free compact
closed  category  with biproducts upon some polycategory \catP.
The earlier results of Propositions~\ref{prop:freeCMon} and
\ref{prop:freeMatrix} described how to freely construct the biproduct as
matrices whose elements are drawn from some category \catC, and this
will provide the core of our proof.

Write $CB\catP$ to denote the free compact closed category with biproducts
generated by a compact polycategory \catP.  We refer to the objects of
\catP, their images under $(\cdot)^*$, and the constants \textbf{0}
and $I$ as the \emph{literals} of $CB\catP$.  According to
Corollary~\ref{cor:star-is-bip-nice}, in any compact closed category with
biproducts, $(\cdot)^*$ commutes with
$\oplus$, and since both the biproduct and tensor structures are
freely generated, the objects of $CB\catP$ are formed from the
literals by repeated application of the functors $(-\otimes -)$ and
$(-\oplus -)$\footnote{See \cite{Soloviev:superpositions:1987} for a
  more general treatment of this.}.  
Any object may therefore be described by such a functor
and a vector of literals.

Let $\otimes_n : CB\catP \times \cdots \times CB\catP \to CB\catP$ be the
$n$-fold tensor; similarly let $\oplus_n$ be the $n$-fold biproduct.
Call $N$ a \emph{normal} functor if it is has the form 
\[
N = \oplus_n (\otimes_{m_1}(-),\ldots,\otimes_{m_n}(-)).
\]
\begin{lemma}\label{lem:complete1}
Let  $G$ be a functor constructed from $(-\otimes -)$ and
$(-\oplus -)$; then $G$ is naturally isomorphic to a normal
functor $N_G$
\end{lemma}
\begin{proof}
  The required isomorphism is constructed from the distributivity
  isomorphisms.
\end{proof}

\noindent 
Hence we have that all arrows in $CB\catP$ have the form
\begin{diagram}
  F\overline{A} & \rTo^f & G\overline{B}\\
  \dTo<\iso && \dTo>\iso\\
  N_F\overline{A} & \rTo_{f'} & N_G\overline{B}
\end{diagram}
and since $f'$ is an arrow between normal functors, it has matrix
elements 
\[
f_{ij} : \otimes_{m_i} A_i \to \otimes_{n_j} B_j
\]
each of which is a (possibly empty) sum of arrows from the freely
constructed compact structure, \CirP.

Hence the free compact closed category  with biproducts is produced
by forming $\mathbf{Matr}( \CirP_\mathbb{N}) $---that is, the free
biproduct category on top of the free compact closed category---and
simply adjoining the distributivity isomorphisms.

This free construction is the final piece of category  theory needed in
this article.  In Section~\ref{sec:gener-proof-nets} we'll introduce a
system of proof-nets that represent this category.

\section{Tensor-Sum Logic}
\label{sec:tensor-sum-logic}

In this section we will introduce the syntax of tensor-sum logic in
a sequent calculus \LTS, and give it a semantics over a suitable category.
Let \catA be a category and denote by $F\catA$ the free compact closed
category with biproducts generated by \catA.  The atomic formulae of
\LTS will be the objects of \catA, and the arrows of \catA will give
its non-logical axioms.  In the next section we will generalise to the
situation where the generators are a polycategory, but that
requires a proof-net presentation.  For now we stick to this simpler
case, since the essence of the connectives can be seen equally
well via a sequent presentation.

\begin{definition}\label{def:formula}
  The \emph{formulae} of \LTS are given by the following grammar:
\[
F ::= \ZERO \;|\; I \;|\; A \;|\; A^* \;|\; F \otimes F \;|\; F \oplus F
\]
where $A \in \OBJ{\catA}$ are called \emph{atoms}.  Given a formula
$F$ we define its \emph{de Morgan dual} $F^*$ by:
\begin{gather*}
  \ZERO^* := \ZERO \\
  I^* := I\\
  A^{**} := A \\
  (F_1 \otimes F_2)^* := F_2^* \otimes F_1^* \\
  (F_1 \oplus F_2)^* :=  F_1^* \oplus  F_2^*\;.
\end{gather*}
An \LTS formula is called \emph{multiplicative} if neither $\ZERO$
nor $\oplus$ occur in it.
\end{definition}

\noindent We use the convention that letters $A$,$B$,$C$ etc, range
over the atoms, while $X$,$Y$,$Z$ etc, range over arbitrary formulae.
We take for granted that all formulae are in de Morgan normal
form---that is, with the negation symbol $(\cdot)^*$ occurring only on
atoms.

\begin{definition}\label{def:sequent}
A \emph{sequent} of \LTS has the form 
\[
\Gamma \vd \Delta ; L
\]
where $\Gamma$ and $\Delta$ are lists of formula, respectively called
the \emph{antecedent} and \emph{succedent} of the sequent, and $L$ is
a tree whose leaves are labelled by \emph{loops} from \catA.  Given
two such trees $L_1$, $L_2$, we write $L_1\cdot L_2$ for the tree
formed by fusing their roots;  we write $L_1 + L_2$ for the
tree whose root has $L_1$ and $L_2$ as its only subtrees.  We don't distinguish
between a loop $l$ in \catA and the tree whose only leaf node is $l$.
\end{definition}

\begin{definition}\label{def:proof}
  An \LTS \emph{proof} is a tree of inferences drawn from the rules shown in
  Figure \ref{fig:rules};  the leaves of the tree must be drawn from
  the  axiom group.  A proof is called \emph{multiplicative} if (1)
  only multiplicative formulae occur in it; and, (2) no rule from
  the additive group occurs.  The reduced sequent calculus consisting
  only of multiplicative proofs we call \LT.
\end{definition}

\begin{figure}[t]
\figureline
\begin{center}
\begin{minipage}{100mm}
\noindent
%
%
\textbf{Axiom Group}: where $f: A \to B$ and $h : A \to A$ are arrows
of \catA.

\[
\begin{array}{ccc}
\RuleT{f}{ ($f$-axiom)}{A \vd B \; ; \emptyset}  
& \qquad \qquad \qquad & 
\RuleT{}{ ($h$-unit)}{ \vd \; ; [ h ] }
\end{array}
\]
\vspace{1.5ex}

\textbf{The Cut:}

\[
\RuleT{\Gamma, X \vd \Delta, X \; ; L}{ (cut)}
{\Gamma \vd \Delta \; ; L  }
\]
\vspace{1.5ex}

\textbf{Multiplicative Group:} $\sigma, \tau$ permutations.

\[
\begin{array}{ccc}
\RuleT{ \Gamma \vd \Delta \; ; L \qquad \Gamma' \vd \Delta' \; ; L'}
{ (mix)}
{\Gamma, \Gamma' \vd \Delta , \Delta' \; ; L\cdot L' } 
& \qquad \qquad \qquad &
\RuleT{\Gamma \vd \Delta \; ; L}
{ (exchange)}
{ \tau(\Gamma) \vd \sigma(\Delta) \; ; L}
\\
\\
\RuleT{ \Gamma, X, Y \vd \Delta \; ; L}
{ ($\otimes$-L)}
{ \Gamma, X\otimes Y \vd \Delta \; ; L}
& \qquad \qquad \qquad &
\RuleT{ \Gamma \vd X, Y, \Delta \; ; L}
{ ($\otimes$-R)}
{ \Gamma \vd X\otimes Y, \Delta \; ; L}
\\
\\
\RuleT{ \Gamma\vd \Delta \; ; L}
{ ($I$-L)}
{ \Gamma,I \vd \Delta \; ; L}
& \qquad \qquad \qquad &
\RuleT{ \Gamma \vd \Delta \; ; L}
{ ($I$-R)}
{ \Gamma \vd \Delta,I \; ; L}
\\
\\
\RuleT{ \Gamma\vd \Delta,X \; ; L}
{ (*-L)}
{  \Gamma,X^* \vd \Delta \; ; L}
& \qquad \qquad \qquad &
\RuleT{ \Gamma, X \vd \Delta \; ; L}
{ (*-R)}
{ \Gamma \vd \Delta,X^* \; ; L}
\end{array}
\]
\vspace{1.5ex}

\textbf{Additive Group}: where $i=1$ or 2.

\[
\hspace{-.75cm}\begin{array}{@{}ccc}
\RuleT{ \Gamma, X_i\vd \Delta \; ; L}
{ ($\oplus_i$-L)}
{ \Gamma, X_1\oplus X_2 \vd \Delta \; ; L}
& \qquad \qquad \qquad &
\RuleT{ \Gamma \vd \Delta , X_i \; ; L}
{ ($\oplus_i$-R)}
{ \Gamma \vd \Delta , X_1\oplus X_2 \; ; L}
\\
\\
\RuleT{0^X_Y}{ (zero)}{X \vd Y \; ; \emptyset}  
&
&
\RuleT{\Gamma \vd \Delta \; ; L \qquad \Gamma \vd \Delta \; ; L'}
{ (sum)}
{\Gamma \vd \Delta \; ; L+L' }
\end{array}
\]
\caption{Inference Rules for \LTS}\label{fig:rules}
\end{minipage}
\end{center}
\figureend
\end{figure}

\begin{figure} 
\topfigrule
\figureline
\begin{center}
\begin{minipage}{100mm}
\textbf{Axiom Group}: where $f: A \to B$ and $h : A \to A$ are arrows
of \catA.

\[
\begin{array}{ccc}
\RuleT{}{ ($f$-axiom)}{f:A \to B}  
&\qquad \qquad \qquad&
\RuleT{}{ ($h$-unit)}{ \text{Tr}^A_{I,I}(h): I \to I }
\end{array}
\]
\vspace{1ex}

\textbf{The Cut:}
\[
\RuleT{\pi: \Gamma \otimes  X \to \Delta \otimes X}{ (cut)}
{\text{Tr}^X_{\Gamma,\Delta}(\pi):\Gamma \to \Delta }
\]
\vspace{1ex}

\textbf{Multiplicative Group:} $\sigma, \tau$ permutations.

\[
\RuleT{\pi: \Gamma \to \Delta}
{ (exchange)}
{ \sigma\circ \pi\circ\tau^{-1}: \tau(\Gamma) \to \sigma(\Delta)}
\]
\[
\RuleT{ \pi:\Gamma \to \Delta  \qquad \pi': \Gamma' \to \Delta'}
{ (mix)}
{(\pi\otimes \pi'):\Gamma \otimes \Gamma' \to \Delta \otimes \Delta' } 
\]

\vspace{0.5ex}
\begin{center}
(No interpretation for tensor or $I$ rules)
\end{center}
\vspace{-0.5ex}

\[
\begin{array}{c}
\RuleT{ \pi: \Gamma \to  \Delta \otimes X }
{ (*-L)}
{ (\id{\Delta} \otimes \epsilon_X ) \circ (\pi\otimes\id{X^*}) :  
 \Gamma \otimes  X^* \to \Delta}
\\
\\
\RuleT{ \pi:X \otimes \Gamma \to \Delta }
{ (*-R)}
{ (\id{X^*}\otimes \pi)\circ (\eta_X \otimes \id{\Gamma}): 
  \Gamma \to X^* \otimes \Delta}
\end{array}
\]

\vspace{1ex}

\textbf{Additive Group}: where $i=1$ or 2.

\[\hspace{-.25cm}
\begin{array}{c}
\RuleT{ \pi:\Gamma\otimes X_i\to \Delta }
{ ($\oplus_i$-L)}
{ \pi \circ (\id{\Gamma}\otimes p_i ): \Gamma\otimes (X_1\oplus X_2) \to \Delta}
\\
\\
\RuleT{ \pi:\Gamma \to \Delta \otimes X_i }
{ ($\oplus_i$-R)}
{ (\id{\Delta}\otimes q_i)\circ \pi: \Gamma \to \Delta \otimes (X_1\oplus X_2)}
\\
\\
\RuleT{\pi:\Gamma \to \Delta\qquad \pi':\Gamma \to \Delta}
{ (sum)}
{\pi+\pi':\Gamma \to \Delta }
\\
\\
\RuleT{}{ (zero)}{0^X_Y : X \to Y} 
\end{array}
\]
\caption{Semantics for rules of \LTS}\label{fig:semantics}
\end{minipage}
\end{center}
\figureend
\end{figure}

\noindent 
One could summarise the rules of \LTS as ``multiplicative-additive
linear logic with self-dual connectives''.  Certainly one can embed
\MALL into \LTS by translating both multiplicative connectives as
$\otimes$ and both additives as $\oplus $ and nothing will go terribly
wrong.  However, since both connectives of \LTS are self-dual, many cuts
which would be forbidden in \MALL are allowed in \LTS, 
and we must introduce some novel rules to deal with this.
It is worthwhile to point out some of the more idiosyncratic rules.

\paragraph{Axiom Rule}
 In the case that \catA is a discrete category then
  the only arrows are identities so we regain the usual $A\vd A$ axioms.
  The restriction of axioms to ground types is for technical
  convenience; identity axioms for every type are constructable, and
  indeed admissible.

\paragraph{Unit Rule} An distinctive feature of compact
closed categories is the presence of loops, so incorporating this rule
allows an exact connection between the syntax and the semantics to be
established.  Perhaps more importantly, the unit rule allows
``circular'' cuts to be eliminated.
\paragraph{Cut Rule}
The cut rule, as shown here, might be better described as
a trace rule.  The more traditional cut rule,
\[
\RuleT{\Gamma\vd\Delta, A \qquad A,\Gamma'\vd\Delta'}{ (cut)}
{ \Gamma, \Gamma'\vd\Delta, \Delta'}
\]
can be defined in \LTS using the mix and exchange rules, viz:
\[
\RuleT{\Gamma\vd \Delta, A \qquad  \Gamma',A \vd\Delta'}{ (mix)}
{\RuleT{\Gamma, \Gamma', A \vd \Delta, A, \Delta' }{(exchange)}
{\RuleT{\Gamma, \Gamma',A\vd \Delta,\Delta',A }{ (cut)}
{\Gamma, \Gamma' \vd \Delta, \Delta'}}}
\]
\paragraph{Mix Rule}
The mix rule (combined with the two rules for tensor) asserts that the
comma is the same on both sides of the sequent, unlike most 
logics.  It allows usual two-premise cut and tensor rules to be
constructed.
\paragraph{Zero Rule}
Without the zero rule certain cuts are impossible to remove.  It has been
noted that the logic of biproducts is inconsistent: every sequent is
provable.  By including the zero axiom we embrace this inconsistency.
A more computational point of view is that every type is inhabited, at
least by the divergent program, or in the quantum setting, the
evolution with zero probability.
\paragraph{Sum Rule}
This rule asserts that each \LTS proof is a (finite) formal weighted
sum of \LTS proofs, with the weights given by the pair $L$ and $L'$.
Otherwise this rule performs a similar role to the mix, allowing the
usual binary rules for additives to be constructed.

The formulae of \LTS are just the objects of \catA hence we shall
not even bother to distinguish them notationally.  To give semantics
for \LTS its remains to translate proofs into arrows of $F\catA$.

\begin{definition}
  Let $\pi$ be an \LTS proof of the sequent 
\[
X_1,\ldots,X_n \vd Y_1,\ldots Y_m ; L\;.
\]
We define its \emph{denotation}, an arrow  
\[
\denote{\pi} : X_1 \otimes \cdots \otimes X_n \to Y_1 \otimes \cdots
\otimes  Y_m
\]
by recursion over the structure of $\pi$ according to the rules shown
in figure \ref{fig:semantics}.
\end{definition}

\begin{theorem}[Cut Elimination]
  For every \LTS proof $\pi$ of the sequent $\Gamma\vd\Delta ; L$ there exists
  a proof $\pi'$ of $\Gamma\vd\Delta ; L$ which contains no occurrence of
  the cut rule, and such that $\denote{\pi} = \denote{\pi'}$.
\end{theorem}

\noindent The proof proceeds in the  standard way so we omit it here.
(The general strategy can be translated from the proof-net version
presented below.)
We remark that preserving the denotation is the non-trivial part;
otherwise the zero rule can be used to give a cut-free proof
immediately.  

We would like the formulae of \LTS to be in exact correspondence with the
objects of the free category  $F\catA$.  Unfortunately we have several
equations between syntactically distinct formulae.  To work around
this blemish we will introduce a special class of formulae.

\begin{definition}\label{defn:reduced-formulae}
  A formula is called \emph{multiplicatively reduced} if it is different to
  $I$ and contains neither $X\otimes I$ nor $I\otimes X$ as a
  subformula,  for any formula $X$.  A sequent is monoidally reduced if
  all its formula are monoidally reduced.

  A formula is called \emph{additively reduced} if it has no
  subformula of the forms $\ZERO\oplus X$, $X\oplus\ZERO$,
  $\ZERO\otimes X$ or $X\otimes \ZERO$.  A sequent is additively
  reduced if all its formulae are.

  A formula or sequent which is both multiplicatively and additively reduced
  is simply called \emph{reduced}.
\end{definition}

\noindent The content of this definition is that the only place that
$I$ may occur in a reduced formula is under the $\oplus$ connective;
the only reduced formula containing $\ZERO$ is $\ZERO$ itself.  

\begin{proposition}
  Every sequent is provably equivalent to a reduced one.
\end{proposition}
\begin{proof}
  We have the following provable equivalences:
  \begin{gather*}
    X \otimes I \equiv X \\
    I \otimes X \equiv X \\    
    X \otimes \ZERO \equiv \ZERO \\
    \ZERO \otimes X \equiv \ZERO \\
    X \oplus \ZERO \equiv X \\
    \ZERO \oplus X \equiv X
  \end{gather*}
  and the denotation of each proof is an isomorphism in $F\catA$.
  The only remaining case is that of a sequent containing the formula
  $I$;  in this case it can be removed by means of a cut, possibly
  after adjoining a new $I$ on the left or right as needed.
\end{proof}

This result means that the \emph{reduced} formulae of \LTS are in 1-1
correspondence with the objects of $F\catA$.  We call a proof reduced
if its conclusion is reduced.  Note that we
cannot restrict to reduced formulae throughout, since they must be
introduced to construct certain formulae, for example $I \oplus I$.
Having dealt with the objects of $F\catA$ we turn out attention to the
arrows.

\begin{theorem}[Completeness]
  Let $f$ be an arrow in $F\catA$;  there exists a cut-free \LTS proof
  $\pi$ such that  $f = \denote{\pi}$.
\end{theorem}

\noindent We again omit the proof since it follows from a more
general result proved below.  However this theorem marks the end of
the line as far as the sequent system is concerned.  Our attempt to
find a proof-theoretic characterisation of $F\catA$ founders on the
usual curse of sequent calculi: the existence of distinct cut-free
forms for the same proof.  \LTS is especially bad in this respect since it
enjoys a great many sound commuting conversions.

\section{Generalised Proof-Nets}
\label{sec:gener-proof-nets}

In this section we define a system of two-sided proof-nets
constructed over the generators of a compact symmetric polycategory
\catP.  The resulting system of proof-nets will be denoted \pnP.
Using these proof-nets we obtain a logical system closer in behaviour
to a term system: every proof-net has a \emph{unique} normal form.  The
availability of such normal forms allows us to make an exact
correspondence between the \pnP and the free compact closed  category
with biproducts.

\subsection{Tensor-Sum Proof-Nets}
\label{sec:tensor-sum-proof}

Graphical notations for monoidal categories have been studied as far
back as the early seventies
\cite{Kelly1972Many-variable-f,Penrose1971Applications-of,JS:1991:GeoTenCal1}
and such 2-dimensional representations provide for  beautifully
simple reasoning in a setting normally awash with coherence
equations.  When proofs are represented graphically, as in proof-nets
for multiplicative linear logic, a further advantage is gained: by
relaxing the allowed shapes of proofs, from trees to graphs, the
artificial sequentiality imposed by the use of sequent proofs is
removed.  Work on \MLL
\cite{Girard:multiplicatives:1987,DanReg:Multiplicatives:1989,BCST1996Natural-Deducti}
extended the graphical tensor notation to the case of two tensors
which enjoy a ``weak'' distribution law\footnote{Also called
  \emph{linear} distribution.} \cite{CocSee:wdc}.  In these settings the
two tensors are essentially similar, and indeed can be made formally
dual.  In the following we study the case a single self-dual tensor so
while much of the work of \cite{BCST1996Natural-Deducti} applies, we can
make some significant simplifications, and are forced into some
complications too, but the we retain a purely graphical language for
the  multiplicative fragment of \LTS, closely related to the diagrams
of \cite{Coecke2005Kindergarten-Qu}.  We note that because \LTS is so
permissive, no correctness criterion, \`{a} la Danos and Regnier
\cite{DanReg:Multiplicatives:1989}, is needed here.

However, into our multiplicative paradise we must admit the additive
connectives and this complicates matters.  We handle the additive
structure using a system of \emph{slices and boxes}.  The notion of
slices in linear proof-nets first appeared in Girard's 
original \cite{Girard:LinLog:1987} but was not entirely satisfactory for
the  unrestricted multiplicative-additive fragment of linear logic;
the  correctness of the proof-structure as a whole could not be
derived from the  correctness of its slices
\cite{Girard:ProofNets:1996}.  A similar notion was later employed in
\cite{HugVGla:AddProofNets:2003} to give a satisfactory notion of MALL
proof-net.  In the more restricted setting of polarised linear logic,
the naive use of slices works very well since the additional
constraint of polarity forces the additive connectives to cohere
nicely \cite{LaurentSlicing-Polariz}.

Just as compact closed categories are degenerate models of the
multiplicative part of linear logic, the biproduct is a degenerate
version of the linear logic's additive connectives.  Happily, this
degeneracy means that slicing will give good 
results, for essentially the opposite reasons to the  polarised case:
we have so many equations that  the  question of correctness becomes
trivial.  Slices and boxes are defined by mutual recursion.

\begin{definition}\label{defn:slice}
  A \emph{\pnP proof-slice} is a finite directed graph
  with edges labelled by \LTS formulae.  The graph is constructed by
  composing the following links, while respecting the labelling on the
  incoming and outgoing edges.
  \begin{description}
  \item[Premise:] No incoming edges; one outgoing edge.  The edge
    is labelled with an arbitrary formula and the link is unlabelled.
  \item[Conclusion:] One incoming edge; no outgoing edges.  The edge is
    labelled with an arbitrary formula and the link is unlabelled.
  \item[Unit:] No incoming edges; two outgoing edges. The first
    outgoing edge is labelled $X^*$, the other, $X$, for some formula
    $X$. The link itself is labelled by $\eta$.
  \item[Counit:] Two incoming edges; no outgoing edges.  Each counit is
    labelled by $\epsilon $ and its incoming edges are labelled by $X$
    and $X^*$ for an arbitrary formula $X$.
  \item[Tensor:] Two incoming edges labelled $X$ and $Y$; one outgoing
    edge labelled $X\otimes Y$.
  \item[Cotensor:]\, One incoming edge labelled $X\otimes Y$; two
    outgoing edges labelled $X$ and $Y$.
  \item[Circle:]\, No incoming or outgoing edges; a circle is a closed
    loop labelled by a formula.
  \item[Axiom:]\, Each polyarrow $f: \langle A_i\rangle_i \to \langle
    B_j \rangle_j$ in
    \ARRC{P} defines a link labelled by $f$.  Its $n$ incoming edges
    are labelled by $A_1,\dots, A_n$ and its $m$ outgoing edges are
    labelled by $B_1,\ldots,B_m$.
  \item[Plus 1:]  One incoming edge labelled $X$; one outgoing edge
    labelled $X\oplus Y$, for an arbitrary formula $Y$.
  \item[Plus 2:]  One incoming edge labelled $Y$; one outgoing edge
    labelled $X\oplus Y$ , for an arbitrary formula $X$.
  \item[CoPlus 1:]  One incoming edge labelled $X\oplus Y$;  one outgoing edge
    labelled $X$.
  \item[CoPlus 2:]  One incoming edge labelled $X\oplus Y$; one outgoing edge
    labelled $Y$.
  \item[Star:]  One outgoing edge labelled $I$;  no incoming edges.
  \item[Costar:]  One incoming edge labelled $I$; no outgoing  edges.
  \item[Box:] Any numbers of incoming and outgoing formulae
    edge, labelled by arbitrary formulae---see definition
    \ref{defn:box} below.
  \end{description}
  \begin{figure}[tbp]
    \[
    \begin{array}{ccccccc}
      \includegraphics{images/smproofnets_3.mps} & \quad &
      \includegraphics{images/smproofnets_4.mps} & \quad &
      \includegraphics{images/smproofnets_5.mps} & \quad &
      \includegraphics{images/smproofnets_6.mps}  \\
      \text{Premise} && \text{Conclusion} &&\text{Unit} &&\text{Counit}   \\ \\
      \includegraphics{images/smproofnets_1.mps} & \quad &
      \includegraphics{images/smproofnets_2.mps} & \quad &
      \includegraphics{images/smproofnets_7.mps} & \quad &
      \includegraphics{images/smproofnets_8.mps}  \\
      \text{Tensor} && \text{Cotensor} && \text{Circle} &&\text{Axiom}  \\ \\
      \includegraphics{images/newproofnets_4.mps} & \quad &
      \includegraphics{images/newproofnets_5.mps} & \quad &
      \includegraphics{images/newproofnets_6.mps} & \quad &
      \includegraphics{images/newproofnets_7.mps}        \\
      \text{Plus 1} && \text{Plus 2} && \text{Coplus 1} &&\text{Coplus
      2}  \\ \\
      \raisebox{3mm}{\includegraphics{images/newproofnets_3.mps}} & \quad &
      \raisebox{3mm}{\includegraphics{images/newproofnets_2.mps}} & \quad &
      \includegraphics{images/newproofnets_8.mps} & \quad &
          \\
      \text{Star} && \text{Costar} && \text{Box} && \\
    \end{array}
    \] 
    \caption{Links for \pnP Proof-nets}
    \label{fig:gen-pn-links}
  \end{figure}

  A  proof-slice is \emph{oriented} such that edges
  enter the node from the top, and
  exit from the bottom.  This implies that any premise, star or unit link is
  above the links they are connected to and, conversely, any
  conclusion, costar, 
  or counit links are below the links they are connected to.

  The order of premises and conclusions is significant, and the
  \emph{type} of a  proof-slice is the pair $(\Gamma,\Delta$) 
  of lists of formulae determined by the premises and conclusions
  respectively.  Usually this will be written as a sequent
  $\Gamma\vd\Delta$.  The empty slice is valid slice, with type
  $\quad\vd\quad$.

  A premise or conclusion link is called \emph{atomic} if 
  the formula labelling it is a literal; a proof-slice is called
  atomic if all its premises and conclusions are atomic.  A slice is called
  \emph{flat} if it contains no boxes
\end{definition}

Proof-slices are permitted to be disconnected or cyclic, when considered
as directed or undirected graphs.  In particular, an edge may leave a
link and return as an input to the same link, although
the labelling on edges will prohibit this for all except axiom links.  If a
proof-slice is directed-acyclic then it is called
\emph{process-like}.

\begin{definition}\label{defn:box}
  A \emph{\pnP box} is a finite multiset of proof-slices, all of the same
  type;  if its component slices are of type $\Gamma\vd\Delta$ then
  the formula of $\Gamma$ are the \emph{inputs} of the box, and those of
  $\Delta$ are its \emph{outputs}.
\end{definition}

A box may be empty; in which case it may have any inputs and outputs.
Indeed, such an empty box is the only normal proof of the formula
$\ZERO$.

Operationally a box may be viewed as local classical
knowledge (or rather, non-determinism) embedded in one part of the
system --- the distribution of addition over composition codes the
transmission of this information.  The details of this distribution,
presumably  mediated by some classical control structure, will not
investigated here, but it seems  an interesting
direction for further exploration.


\begin{definition}\label{def:depth}
  Let $s$ be a proof-slice;  define its \emph{depth} $d(s)$ as 
  \[
  d(s) = \sum_{i=1}^k d(b_i) + k
  \]
  where the $b_i$ range over the boxes occurring in $s$.  Let $b$ be a 
  box containing slices $\{s_i\}_i$;  then define its depth $d(b)$ by
  \[
  d(b) = \sum_i d(s_i)
  \]
\end{definition}
%

\begin{definition}\label{def:proof-net}
  A \emph{\pnP proof-net} is a \pnP box of finite depth.
\end{definition}

According to Definition \ref{def:proof-net} every slice is contained
in a box, which 
is called the \emph{ambient box} for that slice.  Also, if a slice $s$
contains a box, the slices contained within that box are not
considered part of $s$;  that is, from the point of view of their
containing slice boxes are considered totally opaque.  On the other
hand, the phrase ``a slice in proof-net $\pi$'' should be understood
unrestrictedly, as denoting a slice at any level of nesting with the
proof-net structure.  Since proof-nets have finite depth, any box
occurring within a proof-net is itself a valid proof-net;  hence there
is no loss in generality by assuming that ambient box is always the
top-level.

Since \pnP proof-nets are rather unwieldy objects, it is helpful to
introduce a symbolic shorthand for working with them algebraically.
Writing $s_i$ for a sequence of slices, a box containing those slices
is written as a summation, $\sum_i s_i$.  The crucial ingredient is a
``slice with a hole''.  A hole can be thought of as a link with
arbitrary incoming and 
outgoing edges; we write $s\{\;\}$ to represent a slice with a hole.
Such an object is not a valid part of our notation, we introduce it
only in order write such expressions as $s\{t\}$, where $t$ is a
proof-net fragment having the same incoming and outgoing edges as the
hole in $s$ such that the slice produced by replacing the hole in $s$
with the fragment $t$ is a valid slice.  Simply writing the empty
brackets $\{\;\}$ denotes a slice which all hole -- it has no
structure besides its type.  We use letters $\pi, \pi'$ to denote
proof-nets; $\pi\{ \;\}$ should be understood as a proof-net with a
hole in one of its slices.  When we write $\pi\{\;\}\{\;\}$ to denote
a proof-net with two holes, it should always be understood that both
holes are in the same slice.  It will never be necessary to speak of
holes in separate slices.

\begin{example}
This 2 sliced net encodes the distribution of $\otimes$ over $\oplus$.
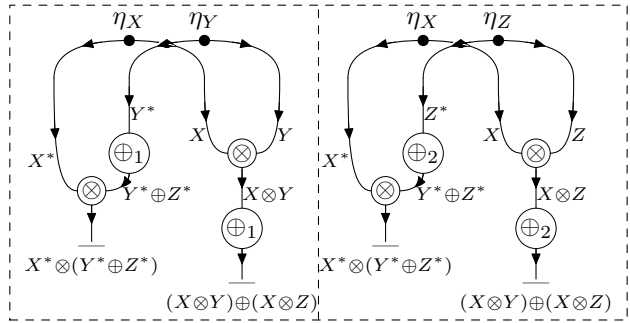
\[
\includegraphics{images/newproofnets_1.mps}
\]
\end{example}

\noindent 
Since categories are a special case of polycategories, we can define
\pnA equally well when \catA is just a category.  In this case the
axiom links have exactly one input and one output;  there is one for
each arrow of \catA.  In this situation we can translate from \LTS
sequent proofs to proof-nets.

\begin{definition}[Translation from sequents]\label{def:trans}
Given an \LTS proof $\pi$, we define a proof-net $N\pi$ by recursion
over the structure of $\pi$.
\begin{itemize}
\item If proof $\pi$ is just an $f$-axiom, let $N\pi$ be the single slice
containing just the corresponding axiom link, connected to a premise and
a conclusion link, leaving a net of type $A\vd B$.
\item If proof $\pi$ is a just an application of the $h$-unit rule for
some $h:A\to A$, we form $N\pi$ by introducing $h$ as an axiom, as
described above, and forming a cut between, as described below,
between its input and output.
\item If $\pi$ is simply an application of the zero rule then $N\pi$
  is an empty box with the desired type.
\item If $\pi$ arises from $\pi'$ by an application of the
cut rule for arrow on some formula $X$ form $N\pi$ from $N\pi'$ by
replacing, in every slice of $N\pi'$, the premise link corresponding
to the negative occurrence of $X$ with an $\eta_X$ link, and
replacing the conclusion link corresponding to the positive occurrence
of $X$ with an $\epsilon_X$ link.  The $X^*$ output of the new unit
link is connected to the $X^*$ input of the new counit link.
\item Suppose $\pi$ arises from subproofs $\pi_1$ and $\pi_2$ by the mix
rule.  Then let $N\pi = \{ N\pi_1\}\{N\pi_2\}$ i.e a single slice
containing two boxes, one for each subproof.
\item If $\pi$ arises from $\pi'$ by an application of the
($\otimes$-R) rule, form $N\pi$ adding, in every slice,  a tensor-link between
the conclusions of $N\pi'$ corresponding to the active formulae of
the rule.
\item If $\pi$ arises from $\pi'$ by an application of the
($\otimes$-L) rule, form $N\pi$ adding, in every slice,  a cotensor-link between
the premises of $N\pi'$ corresponding to the active formulae of
the rule.
\item If $\pi$ arises from $\pi'$ by an application of the
(*-R) rule on some formula $X$, form $N\pi$ adding, in every slice,
an $\eta_X$ -link between to the premise of $N\pi'$ corresponding to
the active formulae of and connect its $X^*$ output to a new conclusion link.
\item If $\pi$ arises from $\pi'$ by an application of the
(*-L) rule on some formula $X$, form $N\pi$ adding, in every slice,
an $\epsilon_X$ -link between to the conclusion of $N\pi'$ corresponding to
the active formulae of and connect its $X^*$ output to a new premise link.
\item If $\pi$ arises from $\pi'$ by an application of the
($I$-R) rule, form $N\pi$ adding, in every slice,
a  star-link, connected a new conclusion link.
\item If $\pi$ arises from $\pi'$ by an application of the
($I$-L) rule, form $N\pi$ adding, in every slice,
a  costar-link, connected a new premise link.
\item If $\pi$ arises from $\pi'$ by an application of the ($\oplus_i$-R)
rule, form $N\pi$ by adding a plus-$i$-link
 to the conclusion corresponding to the active formula in every slice
 of $N\pi'$.
\item If $\pi$ arises from $\pi'$ by an application of the ($\oplus_i$-L)
rule, form $N\pi$ by adding a coplus-$i$-link
 to the premise corresponding to the active formula in every slice
 of $N\pi'$.
\item Suppose $\pi$ arises via an application of the sum rule to
  proofs $\pi_1$ and $\pi_2$; suppose also that $N\pi_1 = \sum_i s_i$
  and $N\pi_2 = \sum_j t_j$.  Then let $N\pi = \sum_i s_i + \sum_j
  t_j$.

\end{itemize}
\end{definition}

\subsection{Normalisation}
\label{sec:normalisation}

\begin{definition}
  Let $e$ be an edge in a slice, going from some link $L_1$ to link
  $L_2$.  We say that  $e$ is \emph{expandable} when: 
  \begin{enumerate}
  \item $e$ is labelled by a a compound formula (\ie either $X\otimes Y$ or
    $X \oplus Y$);
  \item $L_1$ is a premise, cotensor, or coplus link; and, 
  \item $L_2$ is a conclusion, tensor, or plus link.
  \end{enumerate}
\end{definition}

\begin{definition}[Rewrite Steps] \label{def:rewrites}
  Let $\nu,\mu$ be  proof-nets;  define a one step
  reduction relation on proof-nets $R_\beta$ such that  $\nu\,
  R_\beta\, \mu$ if $\nu$ 
  can be rewritten to $\mu$ by one of the following local rewrite
  rules.  
  \subsubsection*{Elimination Rules}\label{rules:elimination}

  \noindent
  \textbf{$\eta\epsilon$-elim:} 
  \[
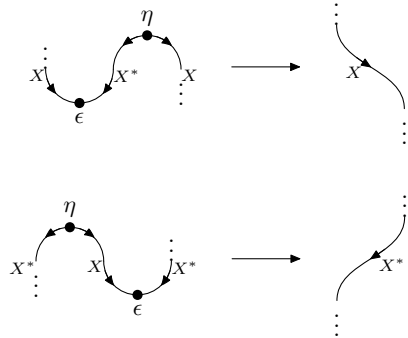
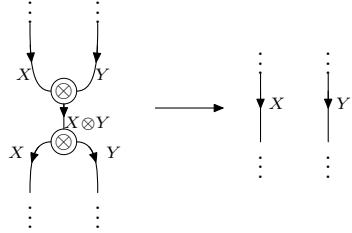
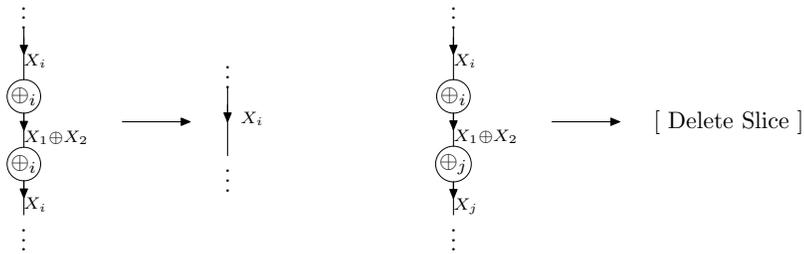
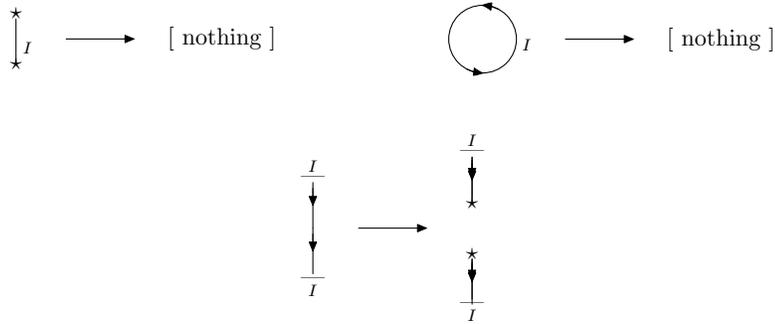
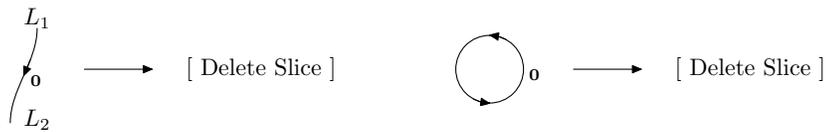
    \begin{array}{c}
    \quad\!\!\raisebox{-0.5\height}{\scalebox{0.9}{
        \includegraphics{images/smproofnets_10.mps}}} \\ \\
    \raisebox{-0.5\height}{\scalebox{0.9}{
        \includegraphics{images/smproofnets_9.mps}}} 
    \end{array}
    \]
    \textbf{$\otimes$-elim}: 
    \[
    \raisebox{-0.5\height}{\scalebox{0.9}{
        \includegraphics{images/smproofnets_11.mps}}} 
    \]

    \textbf{$\oplus$-elim} 
    \[
    \begin{array}{ccc}
      \raisebox{-0.5\height}{\scalebox{0.9}{
          \includegraphics{images/newproofnets_9.mps}}} 
      &\qquad\qquad &         
      \raisebox{-0.5\height}{\scalebox{0.9}{
          \includegraphics{images/newproofnets_10.mps}}} 
    \end{array}
    \]
    where $i \neq j$.

    \noindent    \textbf{$I$-elim}: 
    \begin{gather*}
    \begin{array}{ccc}
    \raisebox{-0.5\height}{\scalebox{0.9}{
        \includegraphics{images/newproofnets_11.mps}}} 
    &\qquad\qquad &             
    \raisebox{-0.5\height}{\scalebox{0.9}{
        \includegraphics{images/newproofnets_12.mps}}} 
    \end{array} \\ \\
        \raisebox{-0.5\height}{\scalebox{0.9}{
        \includegraphics{images/newproofnets_20.mps}}} 
    \end{gather*}
    \noindent \textbf{$\ZERO$-elim} 
    \[
    \raisebox{-0.5\height}{\scalebox{0.9}{
        \includegraphics{images/newproofnets_13.mps}}} \qquad\qquad 
    \raisebox{-0.5\height}{\scalebox{0.9}{
        \includegraphics{images/newproofnets_19.mps}}} 
    \]
    where either of the links $L_1,L_2$ is a
    premise, conclusion, tensor, cotensor, unit, or counit.

    \noindent \textbf{Circle reversal:} 
    \[
    \raisebox{-0.5\height}{\scalebox{0.9}
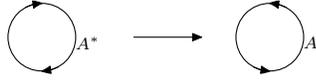
{
        \includegraphics{images/smproofnets_16.mps}}}
    \]
    where $A$ is an atom.

    \subsubsection*{Expansion Rules}\label{rules:expansion}
    \textbf{circle expansion:}
    \[
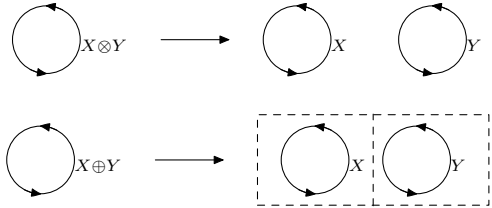
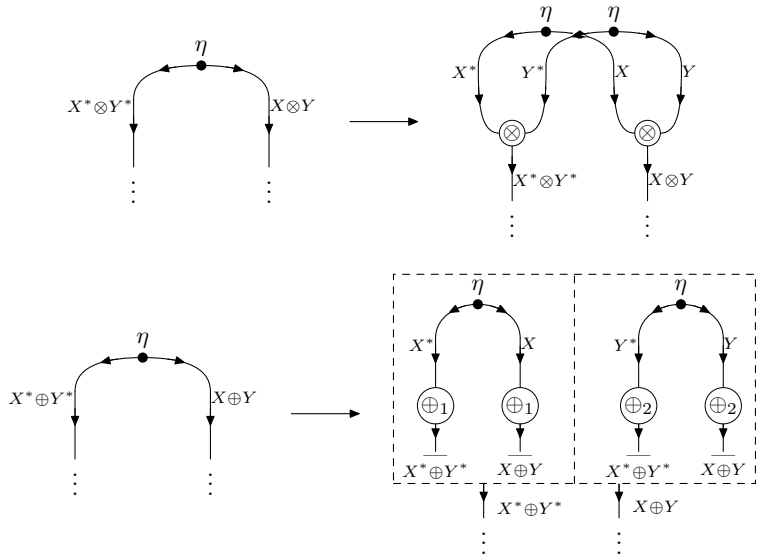
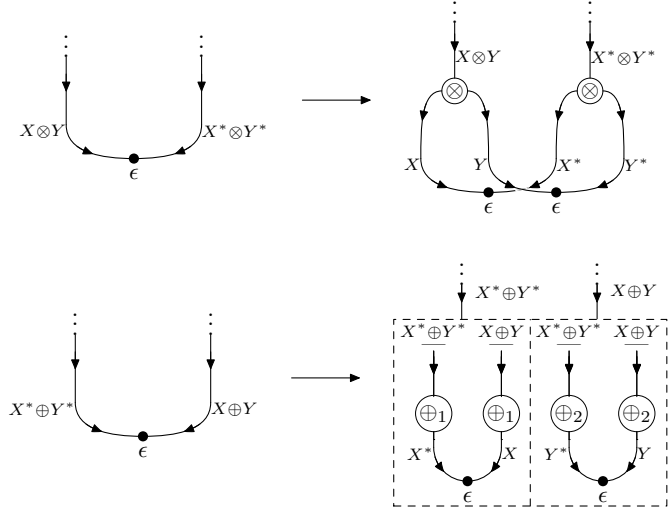
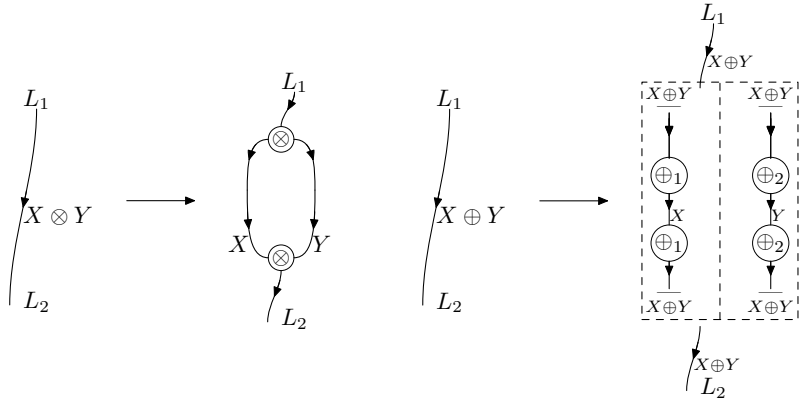
    \begin{array}{c}
    \raisebox{-0.5\height}{\scalebox{0.9}{
        \includegraphics{images/smproofnets_15.mps}}} 
    \\ \qquad \\
    \raisebox{-0.5\height}{\scalebox{0.9}{
        \includegraphics{images/newproofnets_16.mps}}} 
    \end{array}
    \]
    \textbf{$\eta$-expansion:}
    \[
    \begin{array}{c}
    \raisebox{-0.5\height}{\scalebox{0.9}{
        \includegraphics{images/smproofnets_13.mps}}} 
        \\ \qquad \\
    \raisebox{-0.5\height}{\scalebox{0.9}{
        \includegraphics{images/newproofnets_14.mps}}} 
    \end{array}
    \]
    \textbf{$\epsilon$-expansion:}
    \[
    \begin{array}{c}
    \raisebox{-0.5\height}{\scalebox{0.9}{
        \includegraphics{images/smproofnets_14.mps}}} 
        \\ \qquad \\
    \raisebox{-0.5\height}{\scalebox{0.9}{
        \includegraphics{images/newproofnets_15.mps}}} 
    \end{array}
    \]
    \textbf{$\otimes$- and $\oplus$-expansions:} 
    \[
    \begin{array}{ccc}
    \raisebox{-0.5\height}{\scalebox{0.9}{
        \includegraphics{images/newproofnets_17.mps}}} 
    &\quad &             
    \raisebox{-0.5\height}{\scalebox{0.9}{
        \includegraphics{images/newproofnets_18.mps}}} 
    \end{array}
    \]
    An edge from link $L_1$ to link $L_2$ and labelled by a compound
    formula $X$ is expanded when both of the following hold:
    \begin{itemize}
    \item $L_1$ is a cotensor coplus, or premise link; and
    \item $L_2$ is a tensor plus, or conclusion link.
    \end{itemize}

  \subsubsection*{Unboxing Rule}\label{rules:expansion}  
  \noindent If a slice $s$
  contains a box $b = \sum_i t_i$, replace $s$ in the ambient box via
  \[
  s\{\sum_i t_i\} \rTo^\beta \sum_i s\{t_i\}.
  \]
  \ie make a new copy of $s$ for each slice in $b$, and in each
  replace $b$ with the slice.
\end{definition}

\begin{definition}\label{def:beta-reduction}
 Let ${}\rTo_\beta {}$ be the
  transitive, reflexive closure of $R_\beta$ and let  $=_\beta$ be the
  symmetric closure of ${}\rTo_\beta {}$.
\end{definition}

\begin{lemma}[Subject Reduction]\label{lem:gen-pn-subject-reduction}
  Suppose that $\nu$ is a  proof-net with type
  $(\Gamma,\Delta)$ and $\nu \rTo_\beta \mu$; then $\mu$ also has type
  $(\Gamma,\Delta)$.
  \begin{proof}
    No rewrites change the premises or conclusions, hence the type is
    unchanged by $\beta$-reduction.
  \end{proof}
\end{lemma}

\noindent We now begin the approach the proof that $\beta $-reduction is
strongly normalising.  First some intermediary definitions.

\begin{definition}\label{def:expand-depth}
  Let $X$ be a formula; define its \emph{depth} $d(X)$ by
  \begin{gather*}
    d(A) = d(A^*) = d(I) = d(\ZERO) = 1\\
    d(X\otimes Y) = d(X)d(Y) \qquad\qquad d(X \oplus Y) = d(X) + d(Y)
  \end{gather*}
\end{definition}
\begin{lemma}
  Let $e$ be an expandable edge labelled by $X$; then $e$ can be
  expanded to give a box with at most $d(X)$ slices.
\end{lemma}
\begin{proof}
  We use induction on $X$.  Suppose that $X$ contains a connective;
  the expansion rule for that connective will introduce expandable
  edges labelled by the subformulae $Y$ and $Z$.  By induction, these
  yield boxes with $d(Y)$ and $d(Z)$ slices respectively.
  If $X  = Y \otimes Z$; then the expansion rule introduces the new
  edges in parallel; applying the unboxing rule to one
  then the other we, obtain $d(Y)d(Z)$ slices.

  Alternatively suppose $X = Y \oplus Z$.  The expansion rule
  introduces a box containing two slice, with an expandable edge in
  each.  Again, we can apply the unboxing rule twice, and obtain a box
  with $d(Y) + d(Z)$ slices.
\end{proof}

\begin{definition}\label{def:size}
  We define the \emph{size} of a proof-net $\pi$, written $n(\pi)$ by
  mutual recursion over slices and boxes.  Let $s$ be a slice with
  boxes $b_i$;  then
  \[
  n(s) = \prod_i d(b_i) \left( N_s + \sum_i
    \frac{n(b_i)}{d(b_i)}\right)
  \]
  where $N_s$ is the number of links found in $s$, except conclusions,
  premises and boxes.  If we have a box $b = \sum_j t_j$ then let 
  \[
  n(b) = \sum_j n(t_j)\;.
  \]
\end{definition}

\begin{definition}\label{def:rank}
  We define the \emph{rank} of a proof-net $\pi$, written $r(\pi)$ by
  mutual recursion over slices and boxes.  Let $s$ be a slice with
  boxes $b_i$;  then
  \[
  r(s) = \prod_i d(b_i) \left( K_s + \sum_i
    \frac{r(b_i)}{d(b_i)}\right)
  \]
  where $K_s$ is the total number of times the symbols $\otimes$ and
  $\oplus$ occur in the labels of expandable edges of $s$.  If we have
  a box $b = \sum_j t_j$ then let  
  \[
  r(b) = \sum_j r(t_j)\;.
  \]
\end{definition}

\noindent The following lemma is immediate from the definitions:

\begin{lemma}\label{lem:measure-decrease}
  Let $s\{\sum_i t_i\}$ be a slice in a proof-net;  then the following
  hold:
  \begin{gather*}
    n(s\{\sum_i t_i\}) = \sum_i n(s\{t_i\})\\ 
    r(s\{\sum_i t_i\}) = \sum_i r(s\{t_i\}) \\
    d(s\{\sum_i t_i\}) > \sum_i d(s\{t_i\}) \
  \end{gather*}
\end{lemma}

\begin{theorem}[Termination]\label{thm:gen-pn-beta-terminates}
  Every $\beta$-reduction sequence is finite.
\end{theorem}
  \begin{proof}
    We define an order on proof-nets by setting $\pi \succ \pi'$
    whenever $(r(\pi),n(\pi),d(\pi)) > (r(\pi'),n(\pi'),d(\pi'))$ in
    the lexicographic order.  Note that these quantities are all
    non-negative integers so this order has no infinite decreasing chain.

    Suppose now that $\pi R_\beta \pi'$.  
    By inspection of the rules we notice:
    \begin{itemize}
    \item if the rewrite is an expansion, then we have $r(\pi) > r(\pi')$;
    \item if the rewrite is an elimination rule then $n(\pi) >
      n(\pi')$ and $r(\pi) \geq r(\pi')$; and, 
    \item if the  rewrite is the  unboxing rule then by Lemma
      \ref{lem:measure-decrease} we have that $r(\pi) = r(\pi')$,
      $n(\pi) = n(\pi')$ and $d(\pi) > d(\pi')$. 
    \end{itemize}
    Hence $\pi \rTo_\beta \pi'$ then necessarily $\pi \succ \pi'$, and
    therefore every rewrite sequence terminates.
  \end{proof}

\begin{theorem}[Local Confluence]\label{thm:gen-pn-beta-local-confluence}
  If a  proof-net $\pi$ $\beta$-reduces to $\pi_1$ and
  $\pi_2$ by different rewrites $r_1,r_2$, then there exist sequences
  of rewrites $s_1,s_2$ such that 
  \[
  \begin{diagram}
    \pi & \rTo^{r_1} & \pi_1 \\
    \dTo^{r_2} & & \dDashto_{s_1} \\
    \pi_2 &\rDashto_{s_2} & \pi^*
  \end{diagram}
  \]
\end{theorem}

\begin{proof}
  Within a box, rewrites on one slice do not affect any of the other
  slices.  Hence, for a conflict to exist, either $r_1$ and $r_2$ both
  affect the 
  same slice or that  $r_2$ operates on a child slice of that
  where $r_1$ acts. Otherwise there is no conflict between the rewrites and
  they can be trivially unified.  

  The rules also exhibit locality in the vertical
  direction.  The only rule which allows slices on different levels to
  interact is the unboxing rule---and this only pulls slices up from
  the level below.  Hence if $r_1$ and $r_2$ act on slices which are
  more than two levels apart they do not conflict.

  Observe that there are three kinds of rules in the system: those
  that add slices to the ambient box (just the unboxing rule); those
  that delete a slice (the zero rule and the incoherent case of
  $\oplus$-elimination); and those which have purely local effect (all
  the rest).  We'll deal with the cases in that order.
  
  Suppose that the rewrite $r_1$ is the unboxing rule;  without loss
  of generality we have 
  \[
  \pi = \sum_i s_i + s\{\sum_j t_j \} 
  \rTo^{r_1}
  \sum_i s_i + \sum_j s\{ t_j \} = \pi_1
  \]
  Since we need only consider the case where $r_2$ acts on $s\{\;\}$
  or one of the $t_j$, the other slices $s_i$ will be
  neglected. Suppose $r_2$ acts on $s\{\;\}$:
  \begin{itemize}
  \item If $r_2$ is the unboxing rule acting on some other box then we
    have
    \[
    \sum_i s\{t_i\} \{\sum_j t_j'\}
    \lTo^{r_1}
    s\{\sum_i t_i\} \{\sum_j t_j'\}
    \rTo^{r_2}
    \sum_j s\{\sum_i t_i\} \{t_j'\}
    \]
    which can be unified by repeating $r_2$ in each slice on the left,
    and $r_1$ in each slice on the right:
    \[
    \sum_i s\{t_i\} \{\sum_j t_j'\}
    \rDashto^{\sum_j r_2}
    \sum_i \sum_j s\{ t_i\} \{ t_j'\}
    \lDashto^{\sum_i r_1}
    \sum_j s\{\sum_i t_i\} \{t_j'\}
    \]
  \item If $r_2$ deletes $s$ then this must be due some structure in
    $s\{\;\}$ hence the same rule can delete each of the $s\{t_i\}$,
    which suffices to unify the divergence,
  \item If $r_2$ is any other rewrite then we have 
    \[
    \sum_i s\{t_i\} 
    \lTo^{r_1}
    s\{\sum_i t_i\}
    \rTo^{r_2}
    s' \{\sum_i t_i\}
    \]
    where $r_2$ matches some structure in $s\{\;\}$, hence it is still
    available in each of the $s\{t_i\}$, permitting the unification:
    \[
    \sum_i s\{t_i\} 
    \rDashto^{\sum_i r_2}
   \sum_i  s' \{t_i\}
    \lTo^{r_1}
    s' \{\sum_i t_i\}
    \]
  \end{itemize}
  Now suppose $r_s$ acts on one of the $t_i$, which we simply call $t$.
  \begin{itemize}
  \item Suppose $r_2$ is the unboxing rule acting on some box in $t$:
    \[
    t \rTo^{r_2} \sum_j t'_j
    \]
    Then we have the divergence 
    \[
    \sum_{i} s\{t_{i}\} + s\{ t\} 
    \lTo^{r_1}
    s\{\sum_{i} t_{i} + t\}
    \rTo^{r_2}
    s \{\sum_i t_i + \sum_j t'_j \}
    \]
    which we unify using repeated application of the unboxing rule.
    \[
    \sum_{i} s\{t_{i}\} + s\{ t\} 
    \rDashto
    \sum_i s\{t_i\} + \sum_j s\{t'_j\}
    \lTo^
    s \{\sum_i t_i + \sum_j t'_j \}
    \]
  \item Suppose that $r_2$ deletes $t$;  then same rule will delete
    $s\{t\}$, which will unify the divergence.
  \item Otherwise $r_2$ rewrites $t$ to some $t'$;  again this same
    rewrite will do $s\{ t \} \to s\{t'\}$ which will unify the
    divergence.
  \end{itemize}
  This shows that the unboxing rule cannot conflict with the others.
  
  Now suppose that $r_1$ deletes slice $s$.  This implies that $s$
  contains either a pair of incoherent $\oplus$-links or an edge
  labelled by $\ZERO$.  Notice that none of the slice-local rules can
  remove either of these features from the graph.  Hence regardless of
  which rule it is, no slice-local rule $r_2$ can block $r_1$, so the
  divergence can always be unified by deleting $s$.  Of course, If
  $r_2$ also deletes $s$ then there is no divergence.

  Finally we consider the case where both $r_1$ an $r_2$ are
  slice-local. Since all the action is within a single slice, it
  suffices to show that any pair of overlapping rewrites which diverge
  can be unified. 

  Due to the large number of expansion rules, there are a very large number
  of potential critical pairs. Fortunately the rules are very regular and
  have been carefully designed to ensure their confluence.  For
  reasons of space we do not include this analysis here, but checking
  all the pairs is a routine, albeit lengthy, exercise.

  All divergent rewrites can be unified, hence \pnP is locally confluent
  under $\beta$-reduction.
\end{proof}

\begin{theorem}[Strong Normalisation]\label{thm:gen-pn-beta-is-SN}
  $\beta$-reduction for  proof-nets is strongly
  normalising.
\end{theorem}
\begin{proof}
  Since $\beta$-reduction is confluent, each proof-net has a unique
  normal form; since it is terminating, every rewrite sequence must
  arrive at the normal form.
\end{proof}


\noindent
Having established the existence of $\beta$-normal proof-nets, we now
characterise them intrinsically.  Recall that for multiplicative linear logic proof
nets \cite{Girard:multiplicatives:1987},  the structure of a cut free proof
can be separated into the axiom structure and the connective
structure.  The following lemmas give a similar result, 
flattening the additive structure and pushing the
connectives to the outside of the proof-net. 

\begin{lemma}\label{lem:flat-slices}
  Let $\pi$ be a normal proof-net; then every slice of $\pi$ is flat.
\end{lemma}
\begin{proof}
  If any slice contains a box, we can apply the unboxing rule,
  contradicting the normality of $\pi$.
\end{proof}

\noindent Since the box structure of a normal net is trivial, we turn
our attention to the structure of the  slices.  Notice that we will
assume that a normal slice is in a normal net, and hence the rule for
\ZERO-elimination has been applied, implying that a normal slice
contains no edge labelled by \ZERO.

\begin{lemma}\label{lem:gen-pn-notensor}
  Let $\pi$ be a normal proof-slice,  and suppose $x$ is a
  link in $\pi$.  
  \begin{itemize}
  \item If $x$ is a tensor or a plus link, all links below $x$ are
    tensors, pluses, or conclusions.
  \item If $x$ is a cotensor or coplus link, all links above $x$ are
    cotensors, copluses, or premises.
  \end{itemize}
\end{lemma}
  \begin{proof}
    Let $x$ be either a tensor link, or a plus link.  Its
    outgoing edge is labelled by some formula, either $X\otimes Y$ or
    $X \oplus Y$; suppose there is a link below it, called $x'$.  Note
    that since $\pi$ is normal, $x'$ cannot be a box.
    \begin{itemize}
    \item If $x'$ is a counit then it is labelled by a non-atomic
      formula, hence an $\epsilon$-expansion rewrite applies and $\pi$ is
      not normal.
    \item If $x'$ is an axiom, it has an incoming edge labelled by a
      non-atomic formula, which contradicts the definition of axiom link.
    \end{itemize}
    Now there are two cases depending on what kind of link $x$ is.
    \begin{itemize}
    \item Suppose that $x$ is a tensor link;  then $x'$ cannot be a coplus
      link because its incoming formula is $X\otimes Y$.  Suppose that
      $x'$ is a cotensor link: then rewrite rule $\otimes$-elim
      applies, hence $\pi$ is not normal.
    \item Suppose that $x$ is a plus link;  then $x'$ cannot be a cotensor
      link because its incoming formula is $X\oplus Y$.  Suppose that
      $x'$ is a coplus link: then rewrite rule $\oplus$-elim
      applies, hence $\pi$ is not normal.
    \end{itemize}
    Hence $x'$ cannot be a coplus, cotensor, counit, or axiom link.  
    If it is a
    conclusion then the hypothesis is satisfied. If $x'$ is a
    tensor or plus link, then by induction all the links below $x'$ are also
    tensors, pluses, or conclusions.

    The case when $x$ is a cotensor or coplus is exactly dual.
  \end{proof}

\begin{corollary}\label{cor:normal-slice-is-atomic-normal-slice-and-connectives}
  Any normal  proof-slice $\pi$ can be formed from a normal atomic
  slice $\pi'$ by adding tensor and plus links to its conclusions and
  cotensors and copluses  to its premises.
\end{corollary}

\begin{corollary}\label{cor:gen-pn-NF-atomic-labels}
  All the edges of a normal atomic proof-net are labelled by literals.
\end{corollary}

\begin{proposition}\label{prop:gen-pn-atomic-NF-characterisation}
  An atomic proof-slice is normal if and only if: all its edges are
  labelled by atomic formulae; its circles are labelled by positive
  atoms;  no edge is labelled by \ZERO; every edge labelled by $I$
  connects a premise to a costar, or a star to conclusion; and no
  unit link is connected to a counit link;  
\end{proposition}
\begin{proof}
  If $\pi$ is normal, Corollary \ref{cor:gen-pn-NF-atomic-labels}
  gives that all its edges' labels are atomic; by its normality no
  unit is connected to a counit since otherwise rewrite
  $\eta\epsilon$-elim 1 or 2 would apply.  Since $\pi$ is normal, it
  contains no boxes, hence the formula \ZERO can only introduced by a
  premise of conclusion link, but which in this case the \ZERO
  elimination rule would apply.  Since $\pi$ is atomic, the formula
  $I$ may only be introduced by star, costar, premise, or conclusion
  links; any such edge labelled by $I$ is can be eliminated unless it
  connects a premise to a costar, or a star to conclusion as required.

  Conversely, suppose that $\pi$ is atomic, such that  all the above
  conditions are satisfied.
  Since all its edges are
  labelled by literals, none of the expansion rules can apply.  For
  the same reason it contains no tensor, cotensor, plus, coplus, or
  box links, hence 
  rewrites for $\otimes$, $\oplus$, \ZERO, and circle elimination do
  not apply, nor does unboxing.  
  Star and costar links can only appear in forms such that the $I$
  elimination rules do not apply.
  By hypothesis, no unit is connected to a counit, hence rewrites
  $\eta\epsilon$-elim 1 and 2 cannot apply, and circles are labelled by
  positive literals the circle reversal rule does not apply.  Since,
  no rewrites are possible, $\pi$ is in its normal form.
\end{proof}

\subsection{The categorical structure of \pnP}
\label{sec:pna-equivalent-cira}

In this section we prove main remaining theorems about \pnP.  First we
show that \pnP forms a compact closed category with biproducts; and
then we show that it is a representation of the free compact closed
category  with biproducts generated by \catP.

\begin{proposition}
  The class of proof-nets, \pnP, forms a category.
\end{proposition}

\noindent
The objects of \pnP are \LTS formulae.  An arrow $\pi : X\to Y$ is a
proof-net whose only premise is $X$ and whose only conclusion is $Y$.
Two arrows in \pnP are considered equal if they have the same normal
form.

Note that the restriction to single formulae is rather weak since the
comma of tensor-sum logic is implicitly the tensor;
given a proof-net not in this form, we may insert tensor links between
the conclusions, and cotensors between the premises, to obtain a
proof-net of the desired kind.  
The restriction to single formulae also avoids having to provide a
bracketing, since the connectives of \LTS are not strictly
associative.

We define the identity proof-net $\id{X}$ to be a net with one slice,
containing only a premise link and a conclusion link, both labelled by
$X$.  (Note that since the edge linking them may be expandable, this
is not usually the normal form.)

Before defining composition of nets, we first define it for slices.
Suppose $s,t$ are proof-slices such that both the conclusion of $s$, 
and the premise of $t$, are some formula $X$; we define $t\circ s$ by
removing the conclusion link of $s$, removing the premise link of $t$,
and forming a new slice by joining the two graphs along the resulting
open edges.  Notice that this operation is manifestly
associative. Further, we have equations
\[
\id{X} \circ s = s \qquad \text{ and } \qquad t \circ \id{X} = t\;,
\]
since, considering the first case only, we have simply removed a
conclusion link from $s$ and adjoined an identical conclusion link.
The other case is the same.

Now let $f:X\to Y$ and $g: Y\to Z$ be proof-nets, with
slices $f_i$ and $g_j$ respectively;  their
composition $g\circ f = \sum_{ij} g_j\circ f_i$ where
the composition on slices is as above.  Given a third net $h: Z\to W$, 
we have
\[
h\circ(g\circ f) 
= \sum_{ijk} h_k \circ (g_j \circ f_i)
= \sum_{ijk} (h_k \circ g_j) \circ f_i
= (h\circ g)\circ f
\]
so composition of proof-nets is associative as required.  The identity
equations 
\begin{gather*}
\id{Y}\circ f  = \id{Y}\circ (\sum_i f_i) = \sum_i(\id{Y}\circ f_i)  =
\sum_i f_i = f \\
f\circ\id{Y}
 = (\sum_i f_i) \circ \id{Y} 
 = \sum_i (f_i \circ \id{Y} )
 = \sum_i f_i 
= f
\end{gather*}
follow directly from the slice case.  Hence all the axioms required to
be a category are satisfied.

\begin{proposition}
  \pnP is compact closed.
\end{proposition}

\begin{figure}[h]
  \centering
  \figureline\newline
  \[
 \rotatebox{90}{\raisebox{-0.5\height}{\scalebox{0.9}{
          \includegraphics{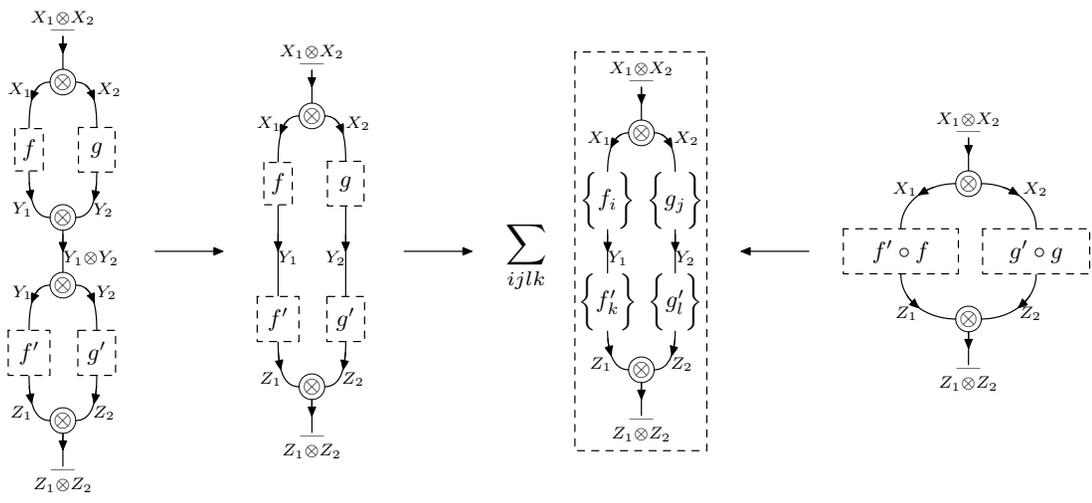}}}} 
\]
\caption{Reduction sequence showing that $(f' \otimes g') \circ (f \otimes g) = (f' \circ f) \otimes (g' \circ g)$}\label{fig:redseq}
\figureend
\end{figure}

\noindent
First we define the monoidal structure of \pnP. Let $f :X_1\to
Y_1$ and $g :X_2\to Y_2$ be proof-nets; then define their tensor
product as
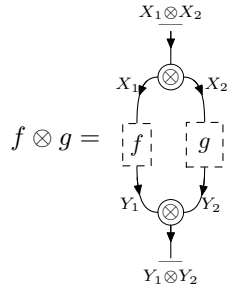
\[
f \otimes g = 
      \raisebox{-0.5\height}{\scalebox{0.9}{
          \includegraphics{images/newproofnets_21.mps}}} 
\]
If $f = \{f_i\}_i$ and $g = \{g_j\}$ then by unboxing we have 
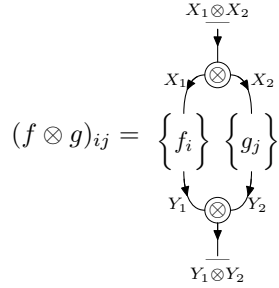
\[
(f \otimes g)_{ij} = 
      \raisebox{-0.5\height}{\scalebox{0.9}{
          \includegraphics{images/newproofnets_22.mps}}} 
\]
Let $f':Y_1\to Z_1$ and $g':Y_2\to Z_2$ be proof-nets then we have the
equation
\[
(f' \otimes g') \circ (f \otimes g) = (f' \circ f) \otimes (g' \circ g)
\]
via the reduction sequence shown in Figure~\ref{fig:redseq}.
To see that $\id{X \otimes Y} = \id{X} \otimes \id{Y}$ we simply
observe that $\id{X \otimes Y} \rTo_\beta \id{X} \otimes \id{Y}$ by
$\otimes$-expansion.  Hence $\otimes$ does indeed define a functor.

The left unit, right unit, symmetry, and associativity isomorphisms are defined
by 
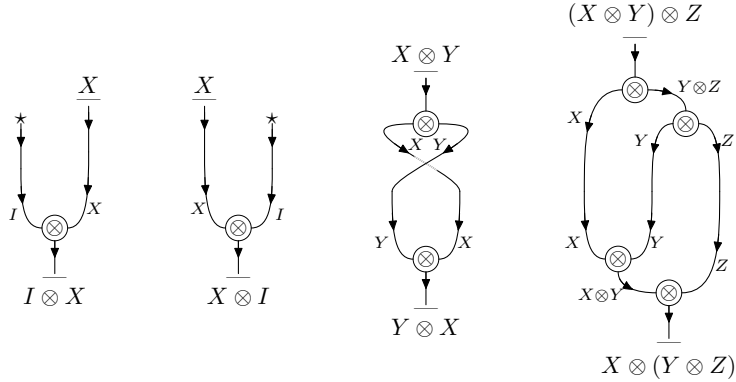
\[
 \raisebox{-0.5\height}{\scalebox{0.9}{
          \includegraphics{images/newproofnets_25.mps}}} 
\qquad \quad
 \raisebox{-0.5\height}{\scalebox{0.9}{
          \includegraphics{images/newproofnets_24.mps}}} 
\qquad \quad
 \raisebox{-0.5\height}{\scalebox{0.9}{
          \includegraphics{images/newproofnets_26.mps}}} 
\qquad \quad
 \raisebox{-0.5\height}{\scalebox{0.9}{
          \includegraphics{images/newproofnets_27.mps}}} 
\]
We leave as an easy exercise to check that the required coherence equations
are satisfied.

Turning attention now to the compact structure, recall that  every
formula $X$ has its de Morgan dual $X^*$ as defined in
Definition~\ref{def:formula}.  The unit and counit maps $\eta_X$ and
$\epsilon_X$ are defined by the nets 
\[
 \raisebox{-0.5\height}{\scalebox{0.9}{
          \includegraphics{images/newproofnets_28.mps}}} 
\qquad \text{ and }\qquad
 \raisebox{-0.5\height}{\scalebox{0.9}{
          \includegraphics{images/newproofnets_29.mps}}} 
\]
The required equations follows more or less immediately from the
$\eta \epsilon$-elimination rules.  Hence $\pnP$ is compact closed.

\begin{proposition}\label{prop:pnP-Cmon-enriched}
  \pnP is enriched over commutative monoids.
\end{proposition}
\begin{proof}
  This property follows directly from the slice structure of
  proof-nets. If $f,g : X\to Y$ are proof nets then $f+g$ is just the
  proof net containing all slices of both $f$ and $g$;  since the
  order of the slices is not significant this operation is
  commutative. The net with no slices, denoted $\emptyset$, gives the
  zero element.
\end{proof}

\begin{proposition}\label{prop:pnP-has-zero}
  \pnP has a \ZERO object.
\end{proposition}
\begin{proof}
  Obviously, the formula \ZERO is the zero object.  Note that for any
  formula $X$, the empty proof-net (i.e the net with no slices)
  provides a proof $\emptyset : X\to \ZERO$ and also $\emptyset :\ZERO
  \to X$.
  
  Suppose that we have a proof-net $f: X\to \ZERO$.  Each slice in $f$
  must  contain a conclusion link labelled by \ZERO; hence by the rule for
  \ZERO-elimination, every slice of $f$ must be deleted, so the normal
  form of $f$ is the empty proof-net.  Hence, for every $X$, there is
  exactly one arrow of type $X \to \ZERO$, and similarly there is
  exactly one arrow $\ZERO\to X$, so \ZERO is both initial and
  terminal in \pnP.
\end{proof}

\begin{proposition}\label{prop:pnP-biproducts}
  \pnP has biproducts.
\end{proposition}
\begin{proof}
  Consider the following one-sliced proof-nets:
  \[
  \pi_1 =  \raisebox{-0.5\height}{\scalebox{0.9}{
      \includegraphics{images/newproofnets_44.mps}}}  
  \qquad 
  \pi_2 =  \raisebox{-0.5\height}{\scalebox{0.9}{
      \includegraphics{images/newproofnets_32.mps}}}  
  \qquad
  \text{in}_1 = \raisebox{-0.5\height}{\scalebox{0.9}{
      \includegraphics{images/newproofnets_33.mps}}}  
  \qquad 
  \text{in}_2 = \raisebox{-0.5\height}{\scalebox{0.9}{
      \includegraphics{images/newproofnets_34.mps}}}  
  \]
  Observe that the rules for $\oplus$-eliminations imply that 
  \[
  \pi_j \circ \text{in}_i = 
  \left\{
  \begin{array}{ll}
    \id{X_i} &\text{ if } i = j \\
    \emptyset &\text{ if } i \neq j
  \end{array} \right. 
  \]
  Next, consider the identity map $\id{X_1\oplus X_2}$.  We have the equation
  \[
  \id{X_1\oplus X_2} = \sum_{i=1,2}   \text{in}_i \circ \pi_i
  \]
  via the rewrite sequence below.
  \[
  \raisebox{-0.5\height}{\scalebox{0.9}
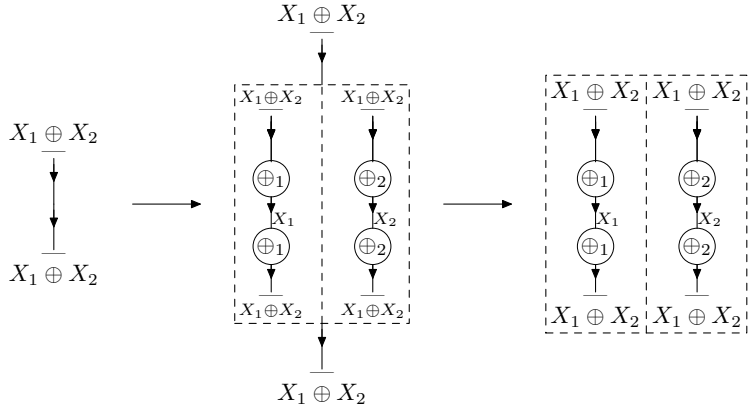
{
      \includegraphics{images/newproofnets_35.mps}}}
  \]
  Since we can form these maps for any pair of objects and, by
  Propositions~\ref{prop:pnP-Cmon-enriched} and
  \ref{prop:pnP-has-zero}, \pnP is a \textbf{CMon}-category with a
  \ZERO object, the result now follows by Proposition
  \ref{prop:CMon-inplies-bips}.
\end{proof}

\begin{proposition}\label{prop:pnP-distribution}
  In \pnP we have natural distribution isomorphisms:
  \begin{gather*}
    X \otimes (Y \oplus  Z) \iso (X\otimes Y) \oplus (X \otimes Z) \\
    (X \oplus Y) \otimes  Z \iso (X\otimes Z) \oplus (Y \otimes Z) \;.
  \end{gather*}
\end{proposition}
\begin{proof}
  The required maps are given by the  proof-nets show below.
  We leave the reader to check that these nets do indeed define
  natural isomorphisms.
  \[
 \raisebox{-0.5\height}{\scalebox{0.9}
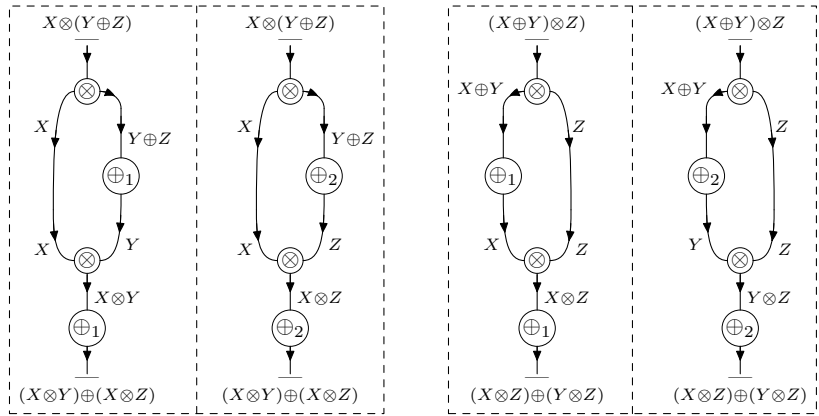
{
          \includegraphics{images/newproofnets_30.mps}}} 
\quad \quad
 \raisebox{-0.5\height}{\scalebox{0.9}{
          \includegraphics{images/newproofnets_31.mps}}} 
  \]
\end{proof}

\noindent
The preceding six propositions established that \pnP is indeed a
compact closed  category  with biproducts as described in Section
\ref{sec:strongly-comp-clos-w-bips}.  Note further that the objects of
\pnP---the \LTS formulae--- are freely generated from the atoms,
which are themselves the objects the underlying polycategory \catP.
Every object of \pnP is therefore isomorphic to a formula in 
disjunctive normal form, 
\[
X \iso \bigoplus_i \otimes_{j_i} X_{j_i}\;,
\]
where the $X_{j_i}$ are literals, and the constants $\ZERO$ and $I$
occur only when a sum or product is empty.  (We assume some given
bracketing of the connectives.)  Hence every proof-net $f : X\to Y$ is
equivalent to some $f'$ of the form:
\[
f' : \bigoplus_i \otimes_{j_i} X_{j_i} 
     \to \bigoplus_{i'} \otimes_{j_i'} Y_{j'_i} \;.
\]
Since $f'$ is a arrow between sums, we can consider its matrix
elements $\pi_i\circ f' \circ \text{in}_j $.  Without loss of generality
take $f'$ to be in normal form;  by Lemmas~\ref{lem:flat-slices} and
\ref{lem:gen-pn-notensor} $f'$ consists of flat slices, whose
connective links are all at the outside, and since its type is in
disjunctive normal form all its plus and coplus links are outside all
its tensor and cotensor links.  Hence 
\[
f' = \sum_{k} \text{in}_{i_k}  \circ f_{k} \circ \pi_{j_k}
\]
where each $f_k$ is a proof-slice between multiplicative formulae. Hence,
\begin{align*}
\pi_i \circ f' \circ \text{in}_j  
& =  \pi_i \circ 
     (\sum_{k} \text{in}_{i_k} \circ f_{k} \circ \pi_{j_k})  
     \circ \text{in}_j \\
& =  \sum_{k}  \pi_i \circ \text{in}_{i_k} \circ f_{k} \circ \pi_{j_k} \circ\text{in}_j\\
& = \sum_{k'} f_{k'}
\end{align*}
where  $k' \in \{ k \; | \; j_k = j \text{ and } i_k = i\}$.  By
Corollary~\ref{cor:normal-slice-is-atomic-normal-slice-and-connectives}
each of the $f_{k'}$ corresponds to a unique normal atomic slice,
which is monoidally reduced.  Hence, the only part the structure of
\pnP which is not freely generated by its connectives are the normal
atomic slices; we now characterise these, and by so doing prove that
\pnP is a representation of the free compact closed  category  with
biproducts generated by \catP.

The reduced normal atomic proof-slices are very closely related to the
\catP-labellable circuits.  Let $\pnP_N$ denote the subcategory of
\pnP determined by the multiplicative formulae, and flat, single-sliced
proof-nets.  We take $\pnP_N$ to be monoidally strict, hence its
arrows are in 1-1 correspondence to the reduced atomic proof-slices.
A simple formal transformation
produces a circuit from each such proof-slice, and vice-versa.  This
correspondence can  be boosted upto a pair of functors
\[
\CirP \pile{\rTo^{F} \\  \\ \lTo_U} \pnP_N
\]
which form an equivalence of categories.

\begin{lemma}\label{lem:gen-pn-NF-negative-edges}
  Suppose $\nu$ is an atomic normal proof-net; suppose $e$ is an edge
  in $\nu$ labelled by a negative literal.  One of the following holds:
  \begin{itemize}
  \item $e$ connects a premise to a conclusion;
  \item $e$ connects a premise to a counit link;
  \item $e$ connects a unit link to a conclusion.
  \end{itemize}
\end{lemma}
\begin{figure}[htbp]
  \centering
  \includegraphics{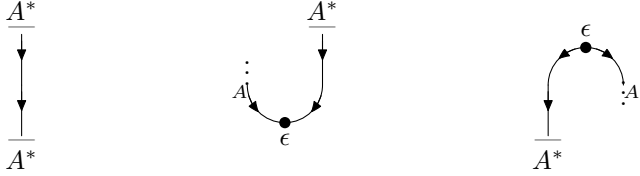}
  \caption{Negative Edges}
  \label{fig:negative-edge-possibilities}
\end{figure}
\begin{proof}
  By Lemma \ref{lem:gen-pn-notensor}, $\nu$ contains no tensor,
  cotensor, plus, or coplus links, nor any boxes; 
  neither axioms nor stars nor costars can introduce negative negative
  edges, therefore $e$ must connect either a premise, unit, counit or
  conclusion.  Since the proof-net is normal, $e$ cannot join a unit
  to a counit by the preceding lemma.  Since an edge cannot be 
  incoming or outgoing at both endpoints the pairings unit/unit,
  counit/counit, premise/premise, conclusion/conclusion,
  conclusion/counit and premise /unit are excluded.  This leaves the
  three possibilities claimed.  These can occur validly in a normal
  proof-net as shown by Fig.~\ref{fig:negative-edge-possibilities}.
\end{proof}

Suppose that $\pi : \Gamma \to \Delta$ is normal and atomic;  then
$\pi$ can be rewritten to produce an \catP-labelled circuit
$c(\pi):\bigotimes\Gamma \to\bigotimes\Delta$
by the following procedure.

\begin{enumerate}
\item The premises and conclusions of $\pi$ become the boundary nodes
  of $c(\pi)$; the premises form $\dom c(\pi)$ and the conclusions
  $\cod c(\pi)$.  They are labelled by the edges formulae and signed
  according to whether the atom is positive or negative.
\item For all edges $e$ labelled by a negative literal $A^*$, 
  reverse $e$'s direction, and change its labelling to $A$.  This
  guarantees that negatively signed nodes in the codomain have
  incoming edges, and vice versa.
\item Erase every unit and counit node, merging their incident edges,
  which are now pointing in the same direction.
\item The remaining links of $\pi$ must all be axioms links.  These
  become the internal nodes of $c(\pi)$.  At each node $x$, the
  ordering on $\text{in}(x)$ and   $\text{out}(x)$ is simply that of
  the components of the domain and  codomain of the arrow (in \catA)
  which labels that node.
\end{enumerate}
Lemma \ref{lem:gen-pn-NF-negative-edges} guarantees that $c(n)$ really
is a circuit.  There is a dual procedure, taking a circuit
$f:\bigotimes_i A_i \to \bigotimes_j B_j$ to a normal atomic proof-net.
\begin{enumerate}
\item The nodes in $\dom f$ become premises; those of $\cod f$,
  conclusions.
\item If $e$ is an edge, labelled by $A$, going from some node $n$ to a
  premise $p$, replace $e$ with a counit-link whose incoming edges are
  from $e$ and $p$, labelled by $A$ and $A^*$ respectively.
\item If $e$ is an edge, labelled by $B$, going to some node $n$ from
  a conclusion $c$, replace $e$ with a unit-link whose outgoing edges
  go to $e$ and $c$ and are labelled by $B$ and $B^*$ respectively.
\item The interior nodes of $f$ become axiom links, each determined by
  the label on the corresponding node.
\end{enumerate}
This defines a proof-net $p(f) : A_1,\ldots,A_n \vd B_1, \ldots B_m$,
which by Proposition~\ref{prop:gen-pn-atomic-NF-characterisation} is
normal.  The two procedures are mutually inverse, which leads to the
following characterisation result.

\begin{figure}[htbp]
  \centering
  \includegraphics[width=\textwidth]{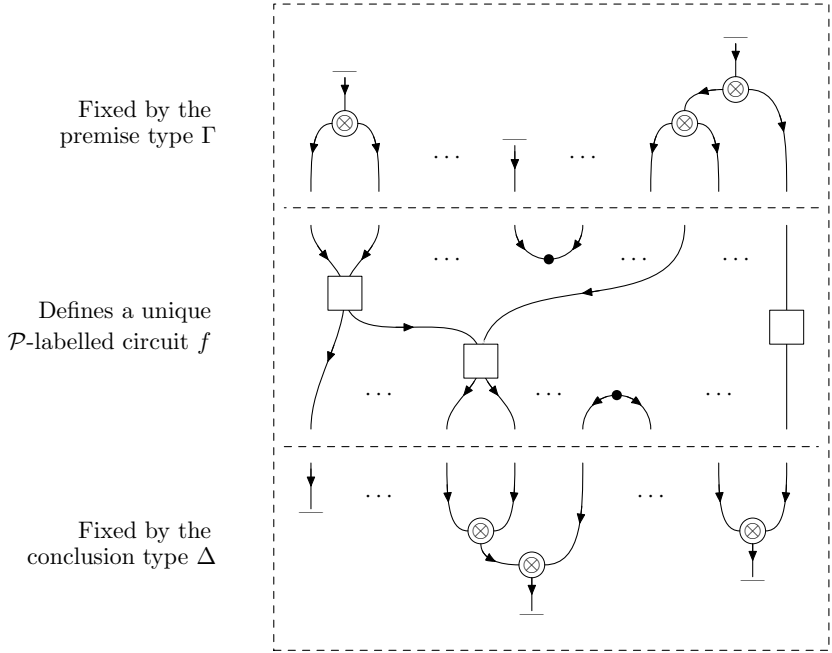}
  \caption{Normal Proof-net decomposition}
  \label{fig:gen-pn-NF-characterisation}
\end{figure}

\begin{definition}
  Let $X$ be formula; an \emph{additive path} for $X$ is a map which
  assigns a boolean value to each occurrence of $\oplus$ in $X$.  
\end{definition}

\noindent
Given an additive path $b$ we can define a purely multiplicative
formula $X_{(b)}$ by replacing each subformula $Y\oplus Z$ with $Y$ if
$b$ assigns 0 to this $\oplus$ and $Z$ if $b$ assigns 1.

\begin{theorem}\label{thm:gen-pn-NF-characterisation}
  Let $\pi$ be a normal proof-slice; then $\pi$ is completely determined
  by its type, an additive path for its domain and codomain, and a
  \catP-labellable circuit.
\end{theorem}
  \begin{proof}
    Suppose that $\pi$ has type $X \vd Y$.  
    Given a formula $X$, and an additive path $b$, let $\langle X,b
    \rangle$ be the list of 
    literals produced by replacing every occurrence of $\otimes$ in
    $X_{(b)}$ by a comma.
    By Lemma \ref{lem:gen-pn-notensor}, $\pi$ can be decomposed into
    three layers:  on top $\pi_{X,b}$ of type $X \vd\langle X,b
    \rangle$
    consisting only of cotensor, coplus, and costar links; the middle $\pi^- :
    X \vd\langle X,b
    \rangle \vd\langle Y,b'
    \rangle$ which is both normal, reduced, and atomic; and at the
    bottom $\pi_{Y,b'}:\langle Y,b' \vd Y$ consisting only
    of tensors, pluses and stars.  The layers  $\pi_\Gamma$ and
    $\pi_\Delta$ are uniquely determined by $X,b$ and $Y,b'$,
    while $\pi^-$ is uniquely determined by the circuit $c(\pi^-)$.
  \end{proof}

  \begin{corollary}
    A normal reduced atomic slice is completely determined by a
    \catP-labelled circuit f.
  \end{corollary}
Justified by the corollary we write $\pi \sim f$ for any normal
reduced atomic slice $\pi$.  The required functors $F$ and $U$ are now
easily defined.

For each $A$ of \pnP, let $UA$ be the positively signed singleton,
labelled by $A$; then define $U(A^*) = (UA)^*$ and $U(X\otimes Y) =
UX\otimes UY$.  Let $\pi \rTo_\beta \nu \sim f$
where $\nu$ is normal; then define $U\pi = f$.

To map \CirP into \pnP, let $f$ be a circuit; then let $Ff$ be the
proof-net obtained from $p(f)$ by adding tensor links to all the
conclusions (bracketed to the  left) and, similarly, cotensors to all
the premises.

\begin{theorem}\label{thm:pnA-CirA-equivalent}
  The 4-tuple $(\CirP,  \pnP_N, F, U)$ is an equivalence of categories.
\end{theorem}
\begin{proof}
  Obviously, from the construction of $U$ and $F$, we have $UF =
  \mathrm{Id}$.  On the other hand, a proof-net $\pi : X\to Y$ only
  differs from $FU\pi: FUX\to FUY$ by the associativity of the
  tensor, hence $\mathrm{Id} \isomorphism  FU$.
\end{proof}

\noindent
This theorem establishes that the matrix elements of a proof-net $\pi$
in \pnP are nothing more than formal sums of circuits over \catP;
i.e. element of the  free compact closed  category  generated by
\catP.  Hence we have the main result:

\begin{theorem}
  The category of proof-nets \pnP is the free compact closed category
  with biproducts generated by the compact symmetric polycategory \catP.
\end{theorem}

\section{Conclusions}
\label{sec:conclusions}

To recap: we sketched how key parts of quantum mechanics can be
formalised in the language of compact closed  categories and
biproducts;  we demonstrated how to represent quantum processes as
proof-nets, and showed that normalisation of such proof-nets allows
some of the behaviour of the corresponding processes to be simulated.

We introduced the formal syntax of tensor-sum logic, and its proof-net
notation.  We showed that proof-nets are strongly normalising, and
characterised the normal forms.  Finally we proved the main theorem: that
the category of proof-nets is exactly the free compact closed  category with
biproducts generated by the polycategory from which its axioms are
drawn.  This result  can be viewed as a coherence theorem for compact
closed  categories  with biproducts, in the style of Kelly and
Laplaza's classic result for compact closed categories
\cite{KelLap:comcl:1980}.  

To return to our starting point, tensor-sum logic is
almost an orthogonal theory to Birkhoff-von Neumann quantum logic.
Tensor-sum logic is entirely preoccupied with the  areas that quantum
logic neglects: compoundness, interaction, and control.  However, as
the main theorem shows,  we
abdicate all responsibility for the internal structure of our quantum
systems.  Since our arrows are characterised by normal proof-nets, 
they are nothing more than type constructors wrapped around the
generators: the fine structure must be described by an equational
theory of the generators.  We can view this as a strength:  the logic is
extremely general and could be easily applied to situations other than
quantum computing.  On the other hand, we suffer strong limitations on
how much of quantum mechanics can be formulated in this setting
without adjoining ad hoc rules to account for the particular
situations we are modelling.

In a sense, this work is the end of the road for those ``logical''
approaches to quantum mechanics deriving from linear logic\footnote{
Those approaches deriving from topos theory, for example
\cite{A.-Doering:2007fu,Heunen:2009qa}, are a different matter
entirely.
}.  Already
the dividing line substructural logic and algebra is thin, and what we
have shown here is that, while proof-theoretic tools may suffice for the
coarse business of putting together systems and pulling them apart
again, the true quantum structure is living in the (poly)category of
generators, and more subtle algebraic tools are needed to tease out
the details.  In particular the importance of spectra in quantum
mechanics weighs against any approach based on natural
transformations.  Recent work
\cite{PavlovicD:MSCS08,Bob-Coecke:2008ty,Coecke:2008nx} provides a
categorical account of observables which is essentially
algebraic.  Fittingly, such theories have graphical representations
which allow them to slot into the proof-net framework as generators.
In that case combining the two systems would yield a well behaved
two-level system of types and terms suitable for representing
quantum processes under classical control.

\bibliography{all}

\end{document}